# Generalized Floquet Exponent, Attractiveness Portrait and Structure Hidden in an Attractor

Keying Guan
Science College, Beijing Jiaotong University,
Beijing, China, 100044
Email: keying.guan@gmail.com

Abstract: The generalized Floquet exponent and the attractiveness portrait (or A-portrait for short) of the attractor and of the smallest invariant closed set are suggested to be used for the study of dynamical systems. Based on the A-portrait, some simple structures hidden in a complicated attractor may emerge from an attractor with complicated structure. The hidden structure plays important role in the bifurcation phenomena of the invariant sets. The examples of A-portraits for the Van der Pol limit cycle, for Lorenz attractor, for the closed limit orbits of different rotation numbers and complicated attractors of Silnikov equation, and for three interlocked smallest invariant closed set of the new improved Nosé-Hoover oscillator are given.

## 1. Some Basic Concepts

In the theory of dynamical systems, the smallest invariant closed set and the attractor are two important conceptions.

Let
$$\frac{dy}{dt} = f(y) \tag{1}$$

be a n-dimensional autonomous differential equation system defined on the open set $D \subset R^n$, where $f$ is a mapping from $D \subset R^n$ to $R^n$ with continuous derivative. The system (1) can be treated as a flow in $D$ with the velocity field $f(y)$.

For a nonlinear system (1), a solution $y(t)$ with initial value $y(0) = \eta \in D$ (denoted by $y(t,\eta)$ ) exists usually in a corresponding time interval $(t^-, t^+)$ where $-\infty \le t^- < 0 < t^+ \le +\infty$. The point sets

$$\gamma = y((t^-, t^+)) \,(= \{y(t) \,|\, t \in (t^-, t^+), y(t) \in D\})$$

$$\gamma^+ = y([0, t^+)) \ (= \{y(t) \mid t \in [0, t^+), y(t) \in D\})$$

and

$$\gamma^- = y((t^-, 0]) \ (= \{y(t) \mid t \in (t^-, 0], y(t) \in D\})$$

are called respectively the orbit (or trajectory), the positive semi-orbit and the negative semi-orbit.

A bounded set $M \subseteq D$ is said to be invariant with respect to the flow (1) if for any solution $y(t)$, its orbit $\gamma$ has an intersection with $M$ (i.e., $\exists t_0 \in R$, $y(t_0) \in M$) implies that $y(t)$ exists for any $t \in R$, and that the whole orbit is in $M$. A nonempty invariant closed set is said to be least or smallest if it has not any nonempty proper invariant closed subset.

A smallest nonempty invariant closed set may be an equilibrium point, a closed orbit, a torus, or a closed set with more complicated structures. This set has usually only a zero measurement, so it has zero possibility (or it is impossible) to find out it numerically if this set has not any attractiveness, or if it has not any stability.

Generally speaking, for a solution $y(t)$, if the positive semi-orbit

$$\gamma^+ = y([0, +\infty)) \ (= \{y(t) \mid t \in [0, +\infty), y(t) \in D\})$$

is located in a bounded closed sub-region of $D$, then the closure $\overline{\gamma^+}$ of this semi-orbit is a non-empty bounded invariant set. If the closed set

$$A(\gamma) = \overline{\gamma^+ \setminus \gamma^+}$$

is nonempty, and if it does not contain $\gamma^+$, then the set $A(\gamma)$ must have some attractiveness to guarantee that $A(\gamma)$ is in the closure of $\gamma^+$. In this case, $A(\gamma)$ is called the attractive set of the orbit $\gamma$.

It is easy to prove that

(i)     The set $A(\gamma)$ is connected,

(ii)    It is a smallest invariant closed set of the system (1)

The detailed proof and the conception "smallest invariant set" can refer to [1] and [2].

For a non-empty smallest invariant closed set $A$, it is called a local attractor of the system (1) if there is an open neighborhood $O(A)$ of it, such that, $y(0) \in O(A)$ implies that the positive semi-orbit $\gamma^+$ of $y(t)$ is in $O(A)$, and that $A$ is the attractive set of $\gamma^+$. This open neighborhood $O(A)$ is usually called an attraction domain of the attractor. Let all of the points in an attraction domain $O(A)$ move along the flow (1), then the domain is getting smaller and smaller. Therefore, the existence of an attractor is related to the phenomenon of the reduction of the volume of the attraction domain when the domain moves along the flow.

For the system (1), the divergence

$$\text{div } f(y)$$

describes the rate of change of the phase volume moving along the flow. Obviously, if the system (1) has an attractor, then there must exist a domain $D_1$ in $D$, such that

$$\text{div } f(y) < 0, \quad \forall y \in D_1$$

and that

$$D_1 \cap A \neq \phi \quad (\phi \text{ is the empty set}) \tag{2}$$

In this case, the system (1) is said to be dissipative in this domain $D_1$.

Since the attraction domain of an attractor has a positive measurement, it is easy to found an attractor numerically by tracing any orbit which starts from a point in its attraction domain.

The system (1) is called conservative if

$$\text{div } f(y) \equiv 0, \quad \forall y \in D$$

A famous conservative system is the 2$n$-dimensional Hamiltonian system

$$\begin{cases} \dfrac{dq}{dt} = \dfrac{\partial H}{\partial p} \\ \dfrac{dp}{dt} = -\dfrac{\partial H}{\partial q} \end{cases}$$

which is generated from the Hamiltonian function $H(q,p)$, $q, p \in R^n$, where $n$ is called the degree of the freedom.

For a conservative system, it cannot have any attractor obviously, but it may have infinitely many smallest invariant closed sets. And a group of smallest invariant closed sets may form a large invariant cluster (set) $\Omega$, which has a positive measurement. For instance, if a Hamiltonian system is integrable by quadratures, or it is a perturbation of a nondegenerated integrable Hamiltonian system, then the system has usually infinitely many invariant $n$-dimensional tori with positive measurement when the perturbation is small. Especially, when the degree of freedom is two, each energy manifold $H=h$ is three-dimensional, the perturbed solutions are either confined on such tori or traped between pairs of such tori (ref. [3], [4]).

Since the measurement of the cluster of these invariant tori is positive, it is easy to find an invariant torus numerically from the cluster. In this sense, the cluster of the invariant tori and the invariant region between a pair of invariant tori in three-dimensional energy manifold is stable.

Just based on the attractiveness of attractor, and the stability of the invariant tori, they can be observed and studied with numerical method.

In mathematics, the Lyapunov exponent (LE) or Lyapunov characteristic exponent of a dynamical system is a quantity that characterizes the rate of separation of infinitesimally close trajectories. So, it can be used to describe the attractiveness and the stability of an attractor or a smallest invariant closed set. Usually, for a *n*-dimensional system, the Lyapunov exponent has *n* components, they also be called the Lyapunov exponents.

By some plausible consideration，for a dissipative system, as criterions, LE is expected to have the following properties (ref. [5]):    if the attractor reduces to

(a) stable fixed point, all the exponents are negative;
(b) limit cycle, an exponent is zero and the remaining ones are all negative;
(c) *k*-dimensional stable torus, the first *k* LEs vanish and the remaining ones are negative;
(d) for strange attractor generated by a chaotic dynamics at least one exponent is positive.

When the phase space is three-dimensional, corresponding to the four different cases (a), (b), (c) and (d) , signs of three Lyapunov exponent components should be distributed respectively as follows:

(a) （-,-,-）
(b) （0,-,-）
(c) （0,0,-）

(d) （+,-,-）, （+,0,-）, （+,+,-）

However, there are still some difficulties and troubles on the Lyapunov exponent in both theory and practical numerical calculation.   One of these difficulties is that there is not a unique exact, objective definition on LE which can realize the above-mentioned properties, (a)-(d). In section 2, some definitions of LE based on the frozen coefficient method are discussed. In section 3, a better definition on LE based on the generalized Floquet exponents is introduced.

The author has noticed that all of the current Lyapunov exponents cannot reflect some profound changes of the structure of an attractor or the smallest invariant closed set when the parameter of the equation changes, such as the period-doubling of a limit closed orbit, though these changes are closely related the attractiveness and stability of the attractor or invariant set. This fact prompts the author to develop some new method to present the attractiveness of these invariant sets in a more detailed way.

In section 4, the attractiveness portrait is introduced for describing the attractiveness distribution on the invariant sets. By this portrait, we may see how and why a limit cycle can be separated into two limit cycle, and how and why a limit closed orbit with a fixed rotation number can change to a limit closed orbit with the doubled rotation number.

In section 5, based on the attractiveness portrait, some subtle structures hidden in certain complicated attractors or smallest invariant sets are revealed. These hidden structures usually related the structures of new invariant sets after the bifurcation of the original complicated invariant set.

## 2. The Lyapunov Exponent Based on the Frozen Coefficient Method

For a given solution $y_0(t)$ of (1), to study its stability or its attractiveness, the system (1) is usually linearized around this solution, i.e.,

$$\frac{dx}{dt} = J(y_0(t))x \qquad (2)$$

where $J(y_0(t)) = \mathrm{D}f\,|_{y=y_0}$ is the Jacobi matrix (or Jacobian for short) of $f$ at $y_0(t)$ (ref. [1]).

If $y_0(t)$ is independent of $t$, i.e., it is the equilibrium solution of (1), then the Jacobian $J(y_0(t))$ is a constant matrix. In this case, Lyapunov gave the strictly criterions to determine the stability of the zero solution with the eigenvalues of the Jacobian. He proved strictly his theory with the method of Lyapunov function.

When $y_0(t)$ depends on $t$ obviously, following the Lyapunov's idea, by freezing the coefficient $J(y_0(t))$ at every fixed $t$, one may study the attractiveness of the zero solution based on the eigenvalues of the Jacobian $J(y_0(t))$ at this moment. By taking different kinds of the average of the local attractiveness, different kinds of Lyapunov exponents can be given. The commonly introduced definitions on LE can be seen in. Wikipedia on "Lyapunov exponent" (http://en.wikipedia.org/wiki/Lyapunov_exponent), or in [5]. They are just based on the frozen coefficient method.

Concretely, for the autonomous system (1), the above-mentioned LEs is to compute first the eigenvales of the average matrix $J(y_0(t))$ for fixed $T\ (>0)$

$$\frac{\int_0^T J(y_0(s))ds}{T} \tag{3}$$

or

$$\frac{\int_0^T J(y_0(s))ds + \text{Transpose}(\int_0^T J(y_0(s))ds)}{2T} \tag{4}$$

then to compute the limits of these eigenvalues as $T \to \infty$.

In [6], the Lyapunov exponents based on (3) was denoted as $LE_J$, the Lyapunov exponents based on (4) was denoted as $LE_O$ (Note: Oseledec proved that $LE_O$ exists with the exception of a subset of initial conditions of zero measure (ref. [5]) ). The present paper will use the same notes.

Recently Zhengling Yang , one of my friend, asked why not to calculate the eigenvalues of $J(y_0(t))$ first before calculating the average (3) or (4). The author cannot find any good reason to response his problem, so in this paper the author suggests another possible definition of LE based on the frozen coefficient method and on Yang's argument:

$$\lim_{T \to \infty} \frac{\int_0^T (\text{Eigenvalues of } J(y_0(s)))ds}{T} \tag{3'}$$

The suggested LE based on (3') will be denoted by $LE_Y$. The numerical results of $LE_Y$ will be also listed in the next section together with $LE_J$ and $LE_O$.

However, there is a serious mathematical problem about these kinds of Lyapunov exponents, that is, if the stability of a solution $y_0(t)$ can be determine completely by the eigenvalues of $J(y_0(t))$ based on the frozen coefficient method.

In fact, this problem has been studied by many authors, such as H.H. Rosenbrook (ref.[7]) , Yuanxun Qin (ref.[8]) and E.W. Kamen et al. (ref.[9]).

In [7], Rosenbrook posed a counterexample,

$$\begin{cases} \dfrac{dx}{dt} = (-1 - 9\cos^2 6t + 12\sin 6t \cos 6t)x + (12\cos^2 6t + 9\sin 6t \cos 6t)y \\ \dfrac{dy}{dt} = (-12\sin^2 6t + 9\sin 6t \cos 6t)x - (1 + 9\sin^2 6t + 12\sin 6t \cos 6t)y \end{cases} \tag{5}$$

Though the Jacobian matrix of (5) has two negative eigenvalues, $\lambda_1 = -1, \lambda_2 = -10$, but the general solution of this system is

$$\begin{cases} x(t) = c_1 e^{2t}(\cos 6t + 2\sin 6t) + c_2 e^{-13t}(\sin 6t - 2\cos 6t) \\ y(t) = c_1 e^{2t}(2\cos 6t - \sin 6t) + c_2 e^{-13t}(2\sin 6t + \cos 6t) \end{cases}$$

Clearly, the zero solution of this system is not stable.

So, generally, the eigenvalues of the Jacobi matrix $J(y_0(t))$ may not globally reflect the stability of the zero solution of (2) when $J(y_0(t))$ is depends obviously on $t$, unless some more constraints are put on (see [7],[8],[9]). Therefore, it should not be expected that the above-mentioned LE$_J$, LE$_O$ or LE$_Y$ could reflect the stability in any case.

In fact, in [6], the author has found that, in the simple case that the attractor is a limit cycle or it is a spatial closed limit orbit, the expectations (a),(b),(c),(d) and (e) on the Lyapunov exponents mentioned in section 1 are not completely correct. The sign distribution of LE$_J$ for limit closed orbit in three dimensional phase space can be any one of the following cases （-,-,-）, （0,-,-）, and （+,-,-）. And in most of cases, the sign distributions of LE$_O$ for a limit closed orbit is （+,-,-）. This means LE$_O$ has almost no sense for the closed orbits.

Therefore, it is in need to find out some better definitions on Lyapunove exponents. In section 3, the generalized local Floquet exponent is suggested.

## 3. Generalized Local Floguet Exponent

In 1883, Achille Marie Gaston Floquet, a French mathematician, developed the famous Floquet theory for the linear periodic differential equation system

$$\frac{dx}{dt} = A(t)x \qquad (6)$$

where the coefficient matrix $A(t)$ is periodic, that is

$$\exists T > 0, \text{ such that } A(t+T) = A(t), \forall t \in R$$

This theory shows that, for any fundamental matrix solution $\phi(t)$ that all columns are linearly independent solutions, there exists a periodic matrix $P(t)$ and a constant matrix $B$, such that

$$\phi(t) = P(t)e^{tB} \qquad (7)$$

If $x \in R^n$, the eigenvalues of the $n \times n$ matrix $B$ are called the Floquet exponents ( ref.[1]

or [10]).

In order to see the meaning of these exponents, choosing a particular case, i.e., these exponents are $n$ different complex numbers, $\mu_1, \mu_2, \ldots, \mu_n$, then it can be proven that there is a corresponding fundamental matrix solution $\hat{\phi}(t)$ of this system, such that,

$$\hat{\phi}(t) = \hat{P}(t) \begin{bmatrix} e^{\mu_1 t} & 0 & \cdots & 0 \\ 0 & e^{\mu_2 t} & \cdots & 0 \\ \vdots & \vdots & \ddots & \vdots \\ 0 & 0 & \cdots & e^{\mu_n t} \end{bmatrix} \quad (8)$$

where $\hat{P}(t)$ is still a periodic matrix. So, the signs of real parts of these exponents determine the stability of the zero solution of the periodic system (6). Sometimes, the real parts of the Floquet exponents is also called the Lyapunov exponents, see

http://en.wikipedia.org/wiki/Floquet_theory

The above-mentioned Floquet exponents reflect exactly the varies of the solutions in a period for the periodic system.

If the linear system (6) is not periodic, the same expression (7) and (8) may exist in any given time interval $[0, T]$. In this case, the aperiodic system (6) can be treated as a periodic system

$$\frac{dx}{dt} = \tilde{A}(t)x, \qquad \tilde{A}(t) = \begin{cases} A(t), & \text{if } 0 \leq t < T \\ A(\tau), & \text{if } t = kT + \tau, 0 \leq \tau < T, k = \pm 1, \pm 2, \ldots \end{cases} \quad (9)$$

The solution of (9) are still continuous, though its derivative may be discontinuous at

$$t = 0, \pm T, \pm 2T, \ldots, \pm kT, \ldots$$

Anyway, the Floquet exponents of the system (7) can exactly reflect the trends of the solutions of the aperiodic system (4) in the time interval $[0, T]$. So, these exponents are called generalized Floquet exponents of the aperiodic system (4) in the local time interval $[0, T]$. This definition can be used for any other time interval $[t_0, t_0 + T]$, $t_0 \in R$.

Generally, it is a difficult problem in the Floquet theory to calculate the exponents exactly with analytical method. However, it is easy to calculate these exponents numerically. One simple way is to calculate the principal fundamental matrix solution $\Phi(t)$, which is a fundamental matrix solution satisfying the initial conditions, $\Phi(t_0) = I$ (the identity matrix). It is easy to prove that the eigenvalues of the matrix $\frac{1}{T} \ln \Phi(t_0 + T)$ are just the generalized Floquet exponents at the local time interval $[t_0, t_0 + T]$.

The generalized Floquet exponent will be denoted by GFE.

In the numerical calculation of the Lyapunov exponents for an attractor of system (1), one must use a numerical solution $\tilde{y}_0(t)$ of (1). It is an elementary knowledge that, the numerical solution $\tilde{y}_0(t)$, no matter how accurate the calculation is performed, it is usually not an exact solution yet, and that, in any given time interval $[t_0, t_0 + T]$, this numerical solution can be treated only as the approximation of a group exact solutions around it.

Therefore, it is almost no sense to calculate numerically the global Lyapunov exponents or generalized Floquet exponents in a very large interval $[t_0, t_0 + T]$, since no one know what the long numerical solution $\tilde{y}_0(t)$ indicates.

The practical numerical calculation for a very long time interval, say $[0, mT]$ ($m$ is a very large positive integer), is to separate this long interval into many ($m$) smaller time intervals with equal length $T$, then to calculate the following averages of all k-th local exponent,

$$\frac{1}{m} \sum_{k=1}^{m} k\text{-th local exponent}$$

respectively.

Concretely, to calculate the average of

all the k-th local eigenvalues of $\dfrac{\int_{(k-1)T}^{kT} J(y_0(s))ds}{T}$ for $LE_J$,

all the k-th local eigenvalues of $\dfrac{\int_{(k-1)T}^{kT} J(y_0(s))ds + \text{Transpose}(\int_{(k-1)T}^{kT} J(y_0(s))ds)}{2T}$ for $LE_O$,

all the k-th local values of $\dfrac{\int_{(k-1)T}^{kT} (\text{Eigenvalues of } J(y_0(s)))ds}{T}$ for $LE_Y$,

and all the k-th eigenvalue of $\frac{1}{T} \ln \Phi_k(kT)$ for GFE.

Note: The calculation of the principal fundamental matrix solution $\Phi_k(kT)$ is based on

$$A(t) = J(\tilde{y}_0(t)), \; t \in [(k-1)T, kT], \; \text{and} \; \Phi_k((k-1)T) = I.$$

It can be proven generally that, for a given n-dimensional system (1) and for a given solution $y_0(t)$,

Summation of $n$ components of LE$_J$ = Summation of $n$ components of LE$_O$
= Summation of $n$ components of LE$_Y$ = Summation of $n$ components of GFE

For all of the summations of $n$ k-th local components of the local LE$_J$, LE$_O$, LE$_Y$ or GFE are equal to

$$\frac{1}{T}\int_{(k-1)T}^{kT} Tr(J(y_0(s)))ds$$

Clearly, the calculation result depends usually on the choice of $T$. It can be proven that LE$_J$, LE$_Y$ and GFE will approach to the same values as $T \to 0$.

One important principle for choosing a better $T$ is naturally that the calculation accuracy for the local exponents can be guaranteed in this time interval.

If the attractor is a closed orbit, the corresponding solution $y_0(t)$ is periodic, these exponents do exist and can be calculated with high accuracy. It is enough to take the calculation in one period. Naturally, it is the best choice to let $T$ be just the period.

Besides the existence of these exponents, the author believes that, the choice for $T$ in the practical calculation is also a serious and principle problem related to the objectiveness of the Lyapunov exponents. It should be studied deeply.

In [11] and [12], the author studied a particular Silnikov equation

$$\begin{cases} \dfrac{dx}{dt} = y \\ \dfrac{dy}{dt} = z \\ \dfrac{dz}{dt} = x^3 - a^2 x - y - bz \end{cases} \quad (10)$$

This system is proven to be an ideal system for the study of three-dimensional differential dynamical systems since it has different kinds of attractors, including spatial limit closed orbits with different rotation numbers.

In [6], as concrete examples, the author calculated a series of LE$_J$ and LE$_O$ for different closed limit closed orbits of a particular Silnikov equation. For different parameters, the numerical results of the four kinds of exponents, LE$_J$, LE$_O$, LE$_Y$ and GFE are listed as follows: (Note: Only the real parts of the corresponding eigenvalues are treated as the exponents )

(n$_1$) $a = 1$ and $b = 0.8$.  The system (10) has a limit cycle of period $T \cong 6.2848$ (Fig.1).

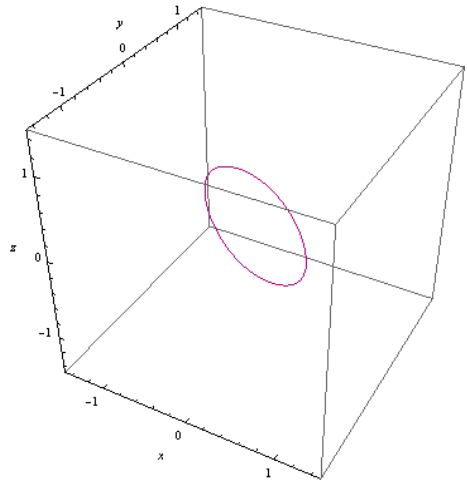

Figure 1. $a = 1$ and $b = 0.8$

For the limit cycle,

| | | |
|---|---|---|
| LE$_J$: | $-0.0701, -0.0701, -0.6600$ | $(-,-,-)$ |
| LE$_O$: | $0.5347, -0.4003, -0.9345$ | $(+,-,-)$ |
| LE$_Y$: | $-0.0935, -0.0935, -0.6130$ | $(-,-,-)$ |
| GFE: | $0.0002, -0,1456, -0.6542$ | $(0^*,-,-)$ |

Note: $0^*$ means that the practical numerical result is approximately equal to zero.

(n$_2$) $a = 1$ and $b = 0.6$.  The system (10) has a limit cycle of period $T \cong 6.2899$ (Fig. 2)

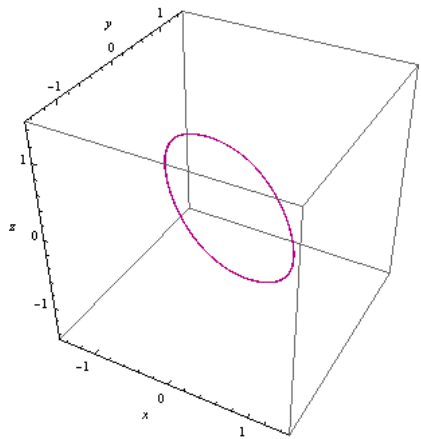

Fig.2   $a = 1$ and $b = 0.6$

For this limit cycle,

| | | |
|---|---|---|
| $LE_J$: | $-0.1929, -0.1929, -0.2142$ | $(-,-,-)$ |
| $LE_O$: | $0.5044, -0.4654, -0.6390$ | $(+,-,-)$ |
| $LE_Y$: | $-0.1958, -0.1958, -0.2085$ | $(-,-,-)$ |
| GFE: | $0.0003, -0.2136, -0.3860$ | $(0^*,-,-)$ |

The following cases (n3)-(n8) seem similar, but there are some slightly changes for the leading (largest) exponents of LEJ and GFE, these changes maybe related the bifurcation of the limit cycle number from one to two. The results will show the GFE seems more exactly sensible to this bifurcation.

($n_3$) $a = 1$ and $b = 0.5024$. The system (10) has a limit cycle of period $T \cong 6.2938$ (Fig. 3).

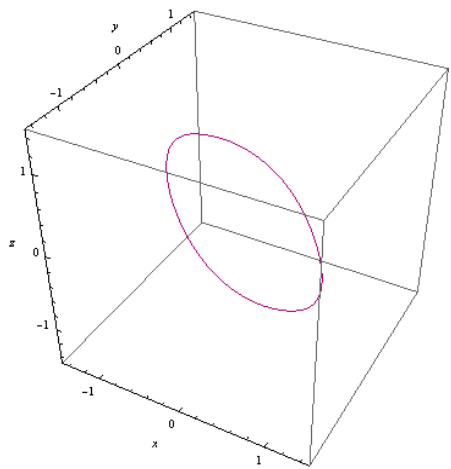

Fig.3    $a = 1$ and $b = 0.5042$

For this limit cycle,

| | | |
|---|---|---|
| $LE_J$: | $-0.00018, -0.2511, -0.2511$ | $(0^*,-,-)$ |
| $LE_O$: | $0.5000, -0.5000, -0.5024$ | $(+,-,-)$ |
| $LE_Y$: | $-0.0692, -0.2166, -0.2166$ | $(-,-,-)$ |
| GFE: | $-0.0022, -0.0207, -0.4794$ | $(0^*,-,-)$ |

($n_4$) $a = 1$ and $b = 0.5023$. The system (10) has still one limit cycle of period $T \cong 6.2938$ (see Fig. 4).   For the limit cycle,

| | | |
|---|---|---|
| $LE_J$: | $0.000017, -0.2512, -0.2512$ | $(0^*,-,-)$ |
| LE$_O$: | 0.5000, $-$0.5000, $-$0.5023 | (+,$-$,$-$) |
| LE$_Y$: | -0.0691, -0.2166, -0.2166 | ($-$,$-$,$-$) |
| GFE: | -0.0022, -0.0206, -0.4795 | ($0^*$,$-$,$-$) |

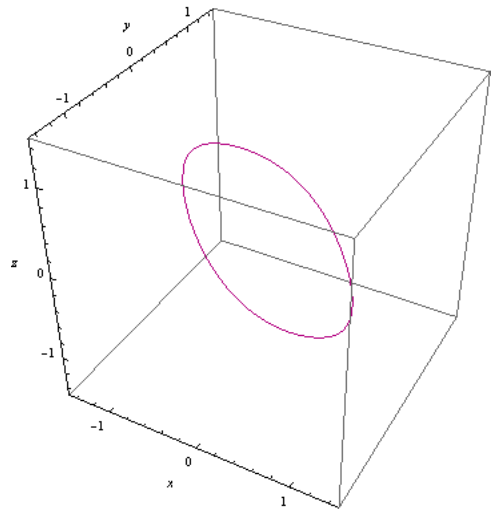

Fig.4  $a = 1$ and $b = 0.5023$

($n_5$) $a = 1$ and $b = 0.5000$. The system (10) has still one limit cycle of period $T \cong 6.2939$ (Fig.5).

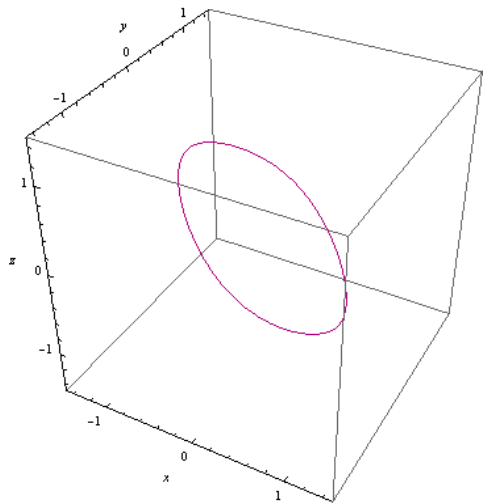

Fig. 5  $a = 1$ and $b = 0.5000$

For this limit cycle,

| | | |
|---|---|---|
| LE$_J$: | 0.0046, $-$0.2523, $-$0.2523 | (+,$-$,$-$) |
| LE$_O$: | 0.5000, $-$0.4984, $-$0.5016 | (+,$-$,$-$) |
| LE$_Y$: | -0.0661, -0.2169, -0.2169 | ($-$,$-$,$-$) |

GFE:    0.0000,  -0.0185,  -0.4815                              $(0^*,-,-)$

($n_6$) $a = 1$ and $b = 0.4900$.  The system (10) has still one limit cycle of period $T \cong 6.2944$ (Fig.6).  For this limit cycle,

LE$_J$:   0.0244, $-0.2572$, $-0.2572$                          $(+,-,-)$
LE$_O$:   0.5000, $-0.4850$, $-0.5051$                          $(+,-,-)$

LE$_Y$:   -0.0532,  -0.2184,  -0.2184                           $(-,-,-)$

GFE:    0.0000,   -0.0013,  -0.4887                             $(0^*,-,-)$

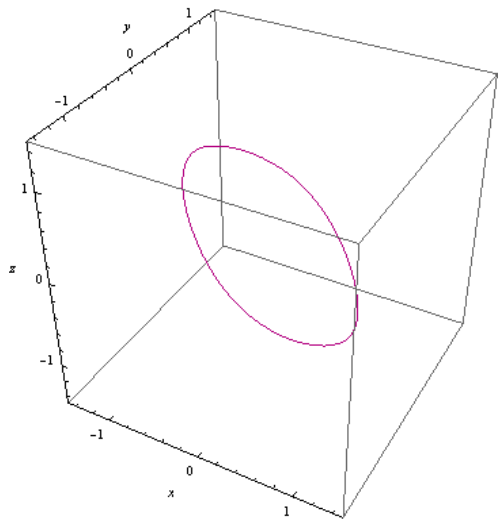

Fig.6  $a = 1$ and $b = 0.4900$

($n_7$) $a = 1$ and $b = 0.4893$.  Numerically, the system (10) has still a limit cycle of period $T \cong 6.2944$ (Fig.7).

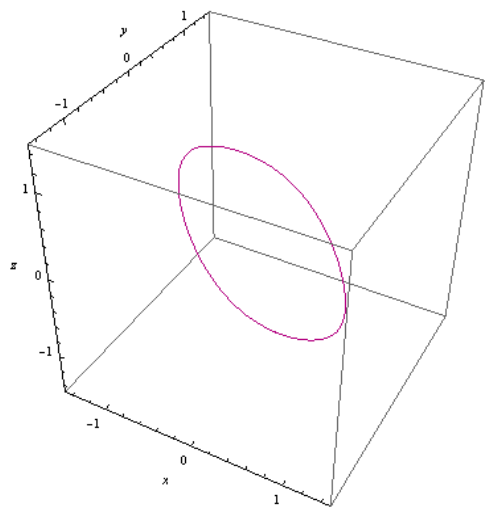

Fig.7   $a = 1$ and $b = 0.4893$.

For this limit cycle,

| | | |
|---|---|---|
| LE$_J$: | 0.0258, −0.2576, −0.2576 | (+,−,−) |
| LE$_O$: | 0.5001, −0.4840, −0.5054 | (+,−,−) |
| LE$_Y$: | -0.0523, -0.2185, -0.2185 | (−,−,−) |
| GFE: | 0.0000, -0.0001, -0.4892 | $(0^*,0^*,−)$ |

(n$_8$)  $a = 1$ and $b = 0.4892$.  From the numerical result, it can be seen carefully that the number of asymptotically stable limit cycles of system (10) has become two (Fig.8). They are symmetrical about the origin and are very close to each other. The period of each limit cycle is $T \cong 6.2944$.  They have the same LE$_J$, the same LE$_O$, the same LE$_Y$ and the same GFE,

| | | |
|---|---|---|
| LE$_J$: | 0.0259, −0.2576, −0.2576 | (+,−,−) |
| LE$_O$: | 0.5001, −0.4839, −0.5054 | (+,−,−) |
| LE$_Y$: | -0.0523, -0.2185, -0.2185 | (−,−,−) |
| GFE: | 0.0000, -0,0002, -0.4875 | $(0^*,0^*,−)$ |

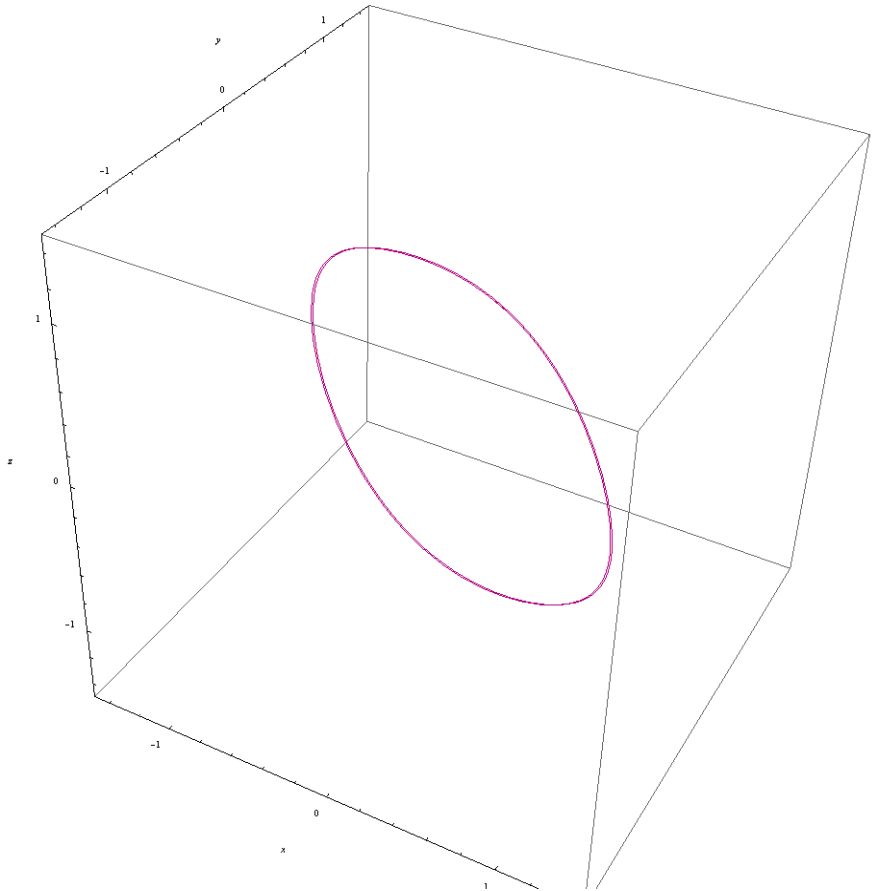

Fig.8  $a = 1$ and $b = 0.4892$,   two slightly separated limit cycles

($n_9$) $a = 1$ and $b = 0.4891$. The system (20) has clearly two asymptotically stable separated and symmetrical limit cycles (see Fig. 9). The period of each limit cycle is $T \cong 6.2945$. For them,

| | | |
|---|---|---|
| LE$_J$:   0.0260, $-0.2576$, $-0.2576$ | | $(+,-,-)$ |
| LE$_O$:   0.5001, $-0.4838$, $-0.5054$ | | $(+,-,-)$ |
| LE$_Y$:  -0.0523,  -0.2184,  -0.2184 | | $(-,-,-)$ |
| GFE:    0.0000,  -0.0006,  -0.4885 | | $(0^*,0^*,-)$ |

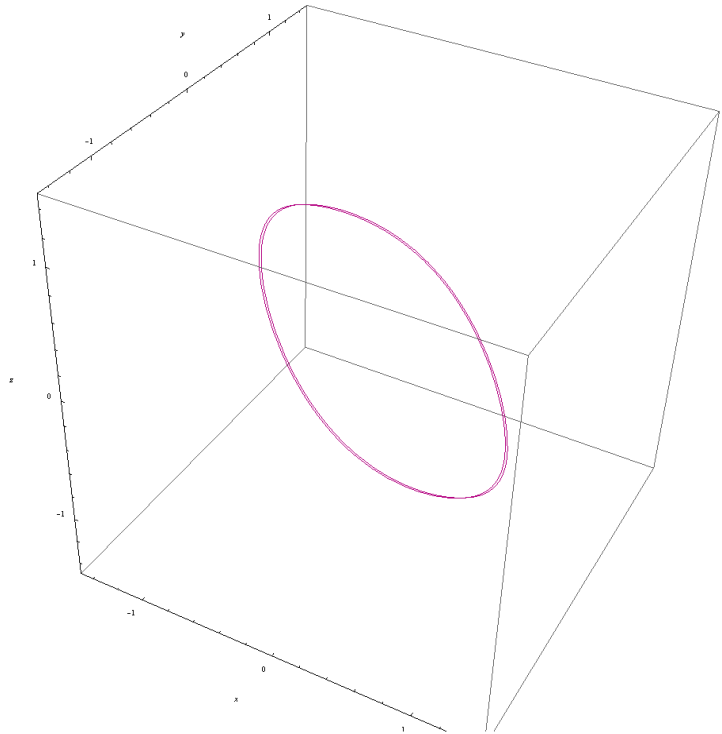

Fig.9  $a = 1$ and $b = 0.4891$, two clearly separate limit cycles

($n_{10}$) $a = 1$ and $b = 0.3920$. The system (10) has two asymptotically stable separated and symmetrical limit closed orbits. Their rotation number are both two (Fig.10). The period of each limit cycle is $T \cong 12.7176$. They have the same LE$_J$, the same LE$_O$, the same LE$_Y$ and the same GFE,

| | | |
|---|---|---|
| LE$_J$:   0.1086, $-0.2503$, $-0.2503$ | | $(+,-,-)$ |
| LE$_O$:   0.5018, $-0.3802$, $-0.5137$ | | $(+,-,-)$ |
| LE$_Y$:  -0.0481,  -0.1720,  -0.1720 | | $(-,-,-)$ |
| GFE:    0.0000,  -0.1960,  -0.1960 | | $(0^*,-,-)$ |

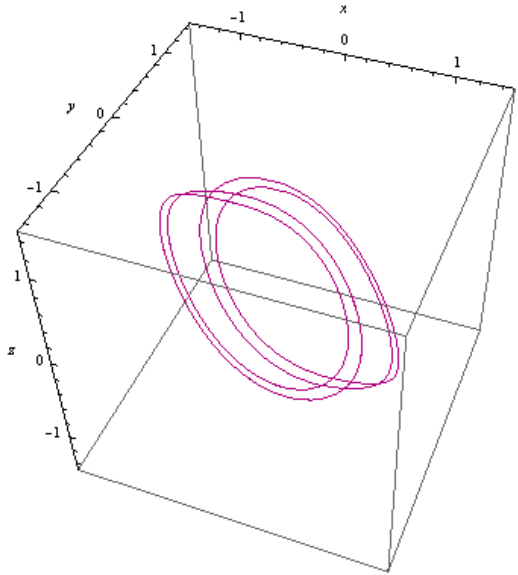

Fig.10  $a = 1$ and $= 0.3920$ , two limit closed orbits of rotation number 2

($n_{11}$) $a = 1$ and $b = 0.3338$.   The system (10) has only one asymptotically stable limit closed orbits. Its rotation number is 13.   It is symmetric to itself about the origin (Fig.11).   The period of the limit closed orbit is $T \cong 84.2365$.

| | | |
|---|---|---|
| $LE_J$:  0.1024,  $-0.2181$,  $-0.2181$ | | $(+,-,-)$ |
| $LE_O$:  0.5031,  $-0.3231$,  $-0.5138$ | | $(+,-,-)$ |
| $LE_Y$:  -0.0580, -0.1379,  -0.1379 | | $(-,-,-)$ |
| GFE:  0.1343,  -0.1016,  -0.3665 | | $(+,-,-)$ |

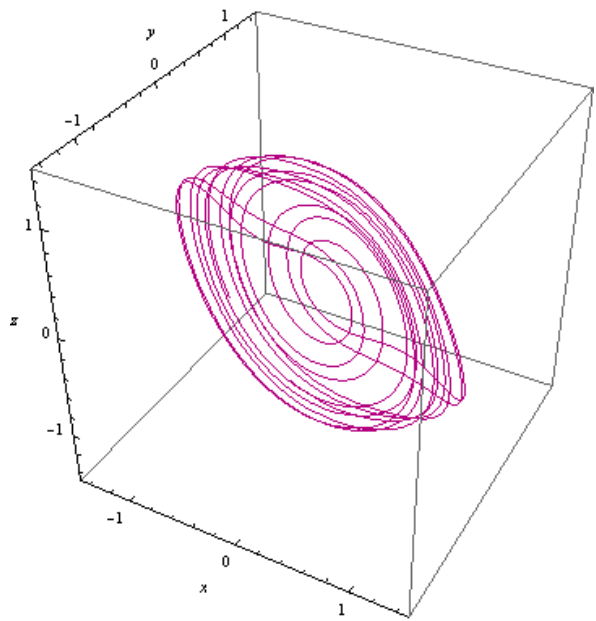

Fig.11  $a = 1$ and $b = 0.3338$. One limit closed orbit with rotation number 13

The above results, ($n_1$)-($n_{10}$), show also that the GFE is more close to the expectation (b) for the Lyapunov exponents listed in section 1 for a limit cycle or a limit closed orbit.

In the case ($n_{11}$), since the period $T \cong 84.2365$ is very large, in order to get a reasonable result, this period has been divided into four smaller intervals in the practical calculation.

The numerical results for the system (10) show that, if the closed orbit has a long period, or it is not a closed orbit, the leading exponent of GFE becomes positive, though the author doubts the accuracy of the calculation for the closed orbits with long period.   See more examples:

($n_{12}$)   When $a = 1$ and $b = 0.3184$, the system (10) has a pair of closed limit orbits.   Both of them have the same rotation number 6 and the same period $T \cong 40.6508$.

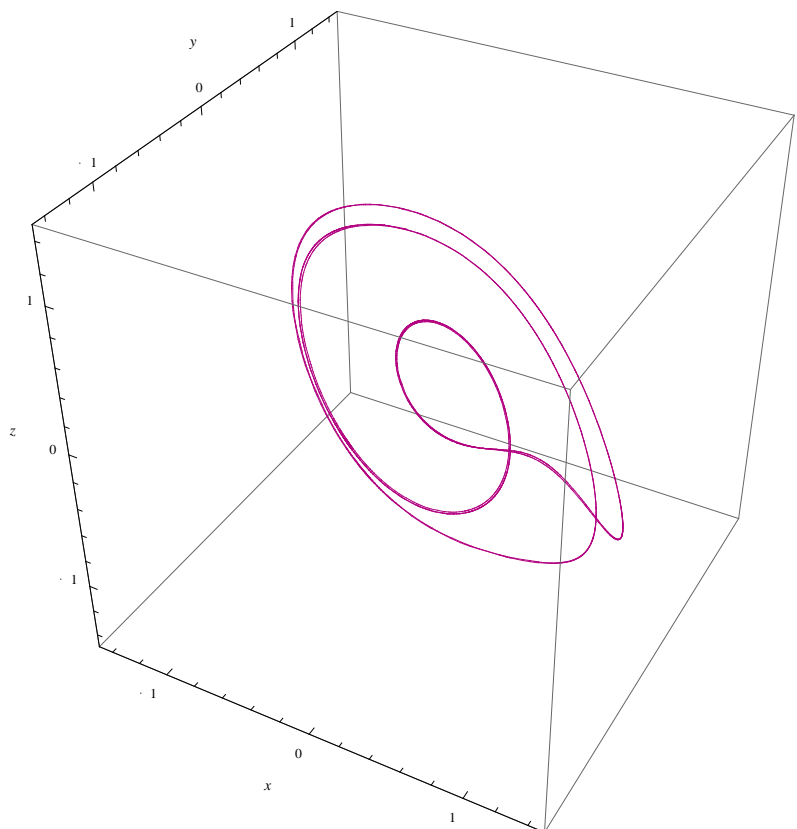

Fig.12  $a = 1$ and $b = 0.3184$, closed orbit of rotation number 6

Figure 12 shows one of the closed orbits with rotation number six.   In this case,

| | | |
|---|---|---|
| $LE_J$:  0.0808,  $-0.1996$,  $-0.1996$ | | $(+,-,-)$ |
| $LE_O$:  0.5011,  $-0.3148$,  $-0.5047$ | | $(+,-,-)$ |
| $LE_Y$:   -0.0887,   -0.1148, -0.1148 | | $(-,-,-)$ |
| GFE:  0.0269,   -0.0313,   -0.3139 | | $(+,-,-)$ |

($n_{13}$)   When $a = 1$ and $b = 0.3140$, the structure of the attractor of the system (10) is complicated, may be just with fractal tructure    (Fig. 13).

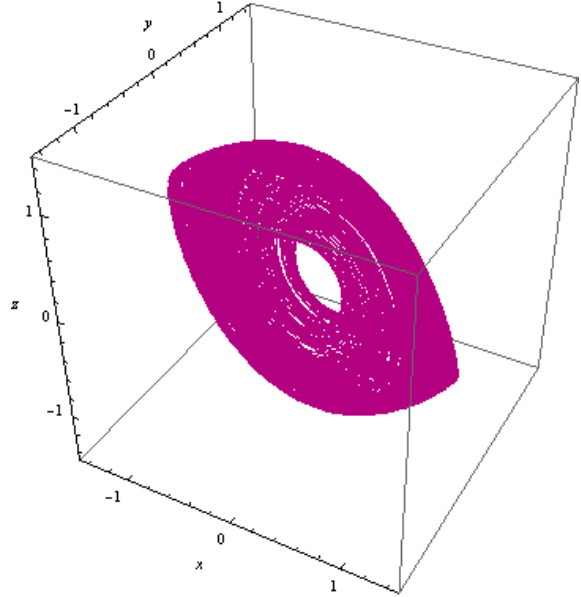

Fig.13   $a = 1$  and  $b = 0.3140$, attractor with complicated structure

There is no period for any trajectory on this attractor.   By choosing $T = 20$, $m = 400$, the numerical results of LE$_J$, LE$_O$, LE$_Y$ and GFE are respectively

| | |
|---|---|
| LE$_J$:   0.1779,  $-0.2460$,  $-0.2460$ | $(+,-,-)$ |
| LE$_O$:   0.5094,  $-0.2910$,  $-0.5325$ | $(+,-,-)$ |
| LE$_Y$:   -0.0045,   -0.1548,  -0.1548 | $(-,-,-)$ |
| GFE:   0.1457,   -0.0378,   -0.4219 | $(+,-,-)$ |

From all of the above examples, it is easy to see that, none of the sign of leading eigenvalue of four different kinds of exponents can satisfy completely the expectation for distinguishing the long period limit closed orbit and attractor with complicated structure.   Besides, it has also been seen that all of these exponents cannot provide enough reliable data for the explanation of the separation of the limit cycle and the explanation of the period-doubling bifurcation of the limit closed orbit.

Therefore, in order to understand the above-mentioned bifurcation phenomena, it is in need to develop some new mathematical methods to present the detailed attractiveness distribution of the attractor or the smallest invariant closed set.   The attractiveness portrait suggested in next section is just an attempt in this direction.

## 4.   Attractiveness Portrait

Rethink the local significance of the idea of the frozen coefficient method. Let the Jacobian $J(y_0(t))$ of the system (1) along a trajectory of the given solution $y_0(t)$ be frozen at a given time $t_0$, then the linearized system (2) can be temporally treated as new autonomous linear system

$$\frac{dx}{dt} = J(y_0(t_0))x \qquad (2')$$

The Lyapunov theory shows that the stability and the attractiveness of the zero solution is determined by the real part of the eigenvalues of the Jacobian $J(y_0(t_0))$. Clearly, this attractiveness should belong also to the solution $y_0(t)$ of the system (1) at $t_0$.

In the following discussion, the number of the dimension of the system (1) is limited as two or three. In the two dimensional case, the eigenvalues of the Jacobian may be two real numbers or a pair of conjugated complex numbers. And in the three-dimensional case, the number of real eigenvalues may be one or three, and if there is only one real eigenvalue, then the other two should be a pair of conjugated complex numbers.

When there is a pair of conjugated complex eigenvalues, there is a corresponding eigenplane on which a trajectory of system (2') is spirally approaching to the origin when the real part of the eigenvalues is negative, or sprirally leaving the origin when the real part of the eigenvalues is positive, or just spirally rotates along a fixed and closed orbit around the origin when the real part of the eigenvalues is zero.

For the real eigenvalues, in most of cases, every real eigenvalue has its own eigenvector, and along the eigenvector, the trajectory of the system (2') is approaching or leaving the origin according to the sign of the eigenvalue. In some particular cases, for instances, the real eigenvalues are douple or triple, the number of the corresponding non-zero eigenvectors may be less than the number of the real eigenvalues, but there is at least one eigenvector corresponding to the double or triple real eigenvalues.

When the real part of the eigenvalue is positive, the zero solution of system (2') exhibits instability or repulsiveness along the corresponding eigenvector or eigenplane. This paper will treat the repulsiveness as a kind of attractiveness, i.e., negative attractiveness. By this understanding, in this paper the word "attractiveness" is usually used to represent both attractiveness and repulsiveness, except in some cases when it is in need to distinguish them.

Based on the above understanding, we may draw an attractiveness portrait (or A-portrait for short) for an attractor or for a smallest invariant closed set of the system (1) according to the following six steps:

$S_1$  Get a numerical solution $y_0(t)$ on the attractor or on the invariant set being studied;

$S_2$  Choose an appropriated positive time length $T$ and such a large positive integer $m$ that the trajectory in the long time interval $[0, mT]$ can be treated as a representation of the attractor or the invariant set.

$S_3$  Calculate the following eigenvalues and the corresponding eigenvectors or eigenplane of $J(y_0(kT))$ for $k = 0,1,2,\ldots,m$

$S_4$  Choose three clearly distinct colors respectively for the trajectory, for the attractive direction and for the repulsive direction. This paper choose blue color to represent attractive for it seems cool, quiescent and convergent, and choose red color to represent repulsive for it seems worm, active, unstable and divergent.  And for the trajectory, choose green.

$S_5$  Draw the trajectory in the phase space with the chosen color (green).

$S_6$  At each point $y_0(kT)$, along every eigenvector of the real eigenvalue, draw a line segment which is centered at the point $y_0(kT)$ with the length proportional to the absolute value of the eigenvalue, and with the attractive color (blue) when the real part of the eigenvalue is negative, or with the repulsive color (red) when the eigenvalue is positive.  If the eigenvalues are a pair of complex numbers, just draw a pair of line segments on the eigenplane, let the two segment cross at their center located at $y_0(kT)$, and let the length of each segment is proportional to the absolute value of the eigenvalues, the line segments are colored with the attractive one (blue) when the real part of the eigenvalue is negative, or colored with the repulsive one (red) when the real part of the eigenvalue is positive.

After the above 6 steps, an A-portrait is obtained. This portrait reflects the distribution and directions of the attractiveness of the attractor or of the invariant set.

Now consider some examples of A-portraits :

$E_1$  A-portrait of Lorenz attractor.  It is well known that the Lorenz equation

$$\begin{cases} \dfrac{dx}{dt} = \sigma(y - x) \\ \dfrac{dy}{dt} = \rho x - y - xz \\ \dfrac{dz}{dt} = xy - \beta z \end{cases} \quad (11)$$

has a strange attractor when $\sigma = 10$, $\beta = \dfrac{8}{3}$ and $\rho = 28$. Figures 14.1, 14.2 and 14.3 show the A-portrait for this attractor from three different viewing points. This portrait is obtained under the choice of $T = 0.4$ and $m = 5000$.

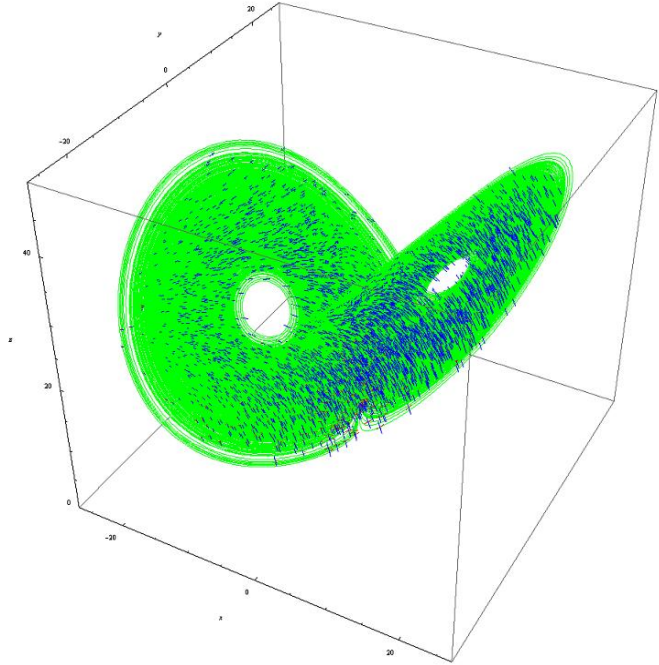

Fig.14.1　A-portrait of Lorenz attractor (1)

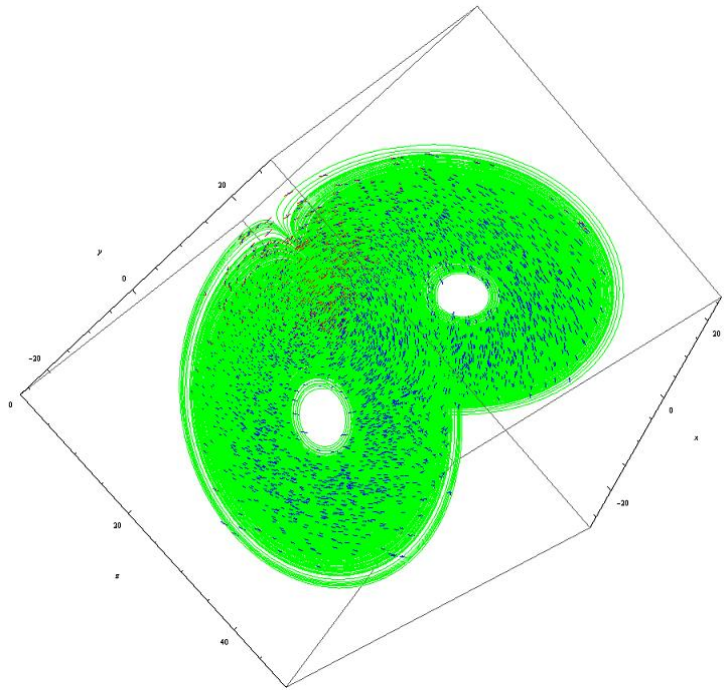

Fig. 14.2　A-portrait of Lorenz attractor (2)

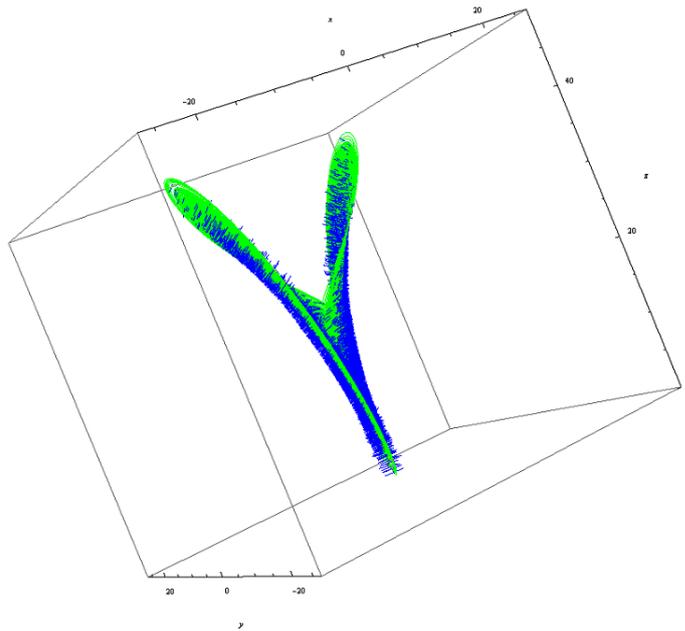

Fig. 14.3    A-portrait of Lorenz attractor (3)

From the A-portrait of Lorenz attractor, just as expected commonly, it can be seen that the attractive directions are mainly perpendicular to the pages of the attractor (like a Cantor book), and the repulsive directions are mainly tangent to the pages of the attractor. The Cantor book is usually used to describe the fractal structure of the strange attractor in three-dimensional phase space (see [4]).    It can also be seen that the distribution of the attractiveness and repulsiveness is uneven on the attractor.

**E$_2$**  A-portrait of a complicated attractor of Silnikov equation.  As shown in Figure13, the Silnikov equation (10) has an attractor with complicated structure when $a=1, b=0.314$.

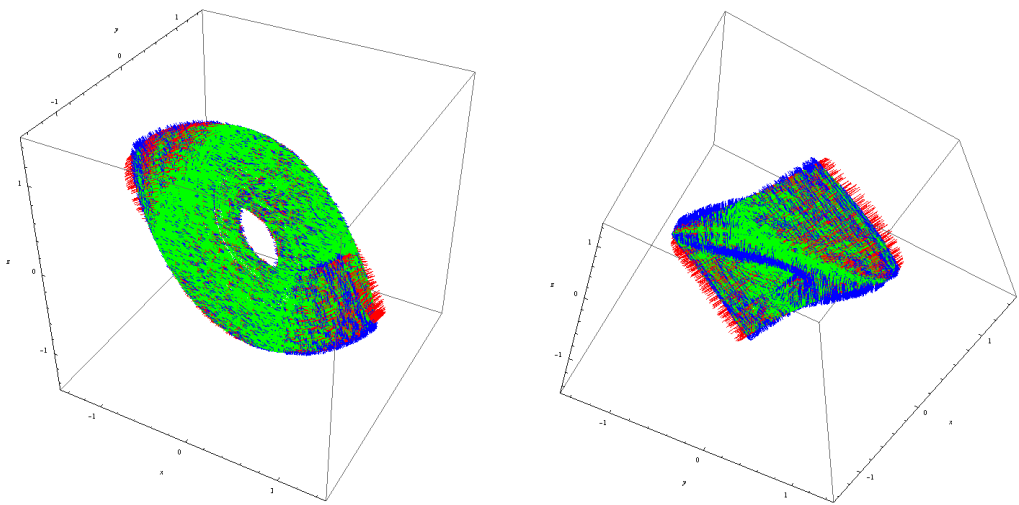

Fig. 15.1    A-portrait of an attractor of the Silnikov equation

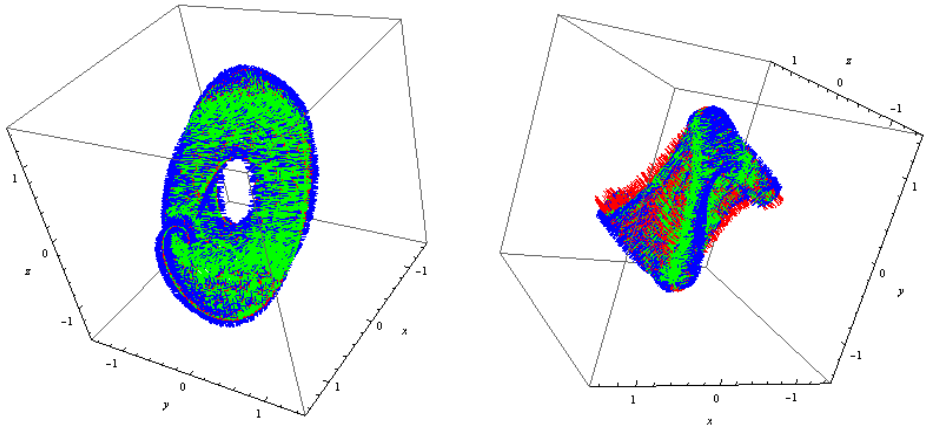

Fig. 15.2   A-portrait of an attractor of the Silnikov equation

Figures 15.1 and 15.2 show the A-portrait from four different viewing points. This portrait is obtained under the choice of $T=1$ and $m=8000$. Similarly to the Lorenz attractor, the attractive directions are mainly perpendicular to the pages of the attractor, and the repulsive directions are just tangent to the pages. Also, the distribution of the attractiveness and repulsiveness is uneven on the attractor.   Note: the obvious repulsive red segments appeared at the edge of the attractor, they are just tangent to the pages of this attractor.

**E₃**   A-portrait of a plane limit cycle.   The following plane system

$$\begin{cases} \dfrac{dx}{dt} = y + x(1-x^2-y^2) \\ \dfrac{dy}{dt} = -x + y(1-x^2-y^2) \end{cases} \quad (12)$$

has a stable limit cycle:  $x=\sin t, y=\cos t$ .   Figure 16 is the A-portrait of the limit cycle.

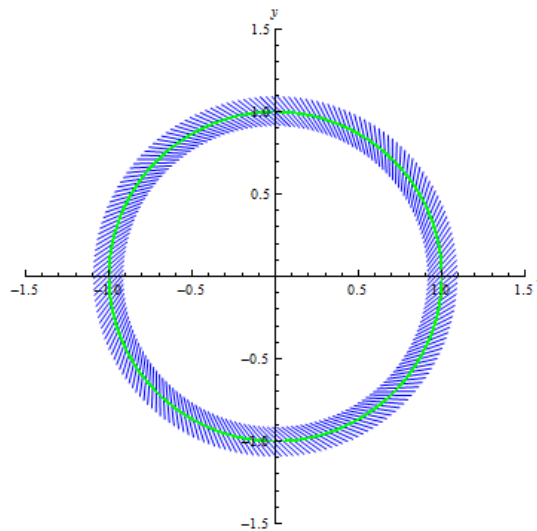

Fig. 16   A-portrait of the limit cycle of (12

Because of the strong symmetry, this limit cycle has only attractiveness distributed evenly.

**E$_4$**    A-Portrait of the limit cycle of the Van der Pol equation

$$\begin{cases} \dfrac{dx}{dt} = y \\ \dfrac{dy}{dt} = -x + y(1-x^2) \end{cases} \quad (13)$$

Figure 17 shows the A-portrait of the limit cycle of this equation.

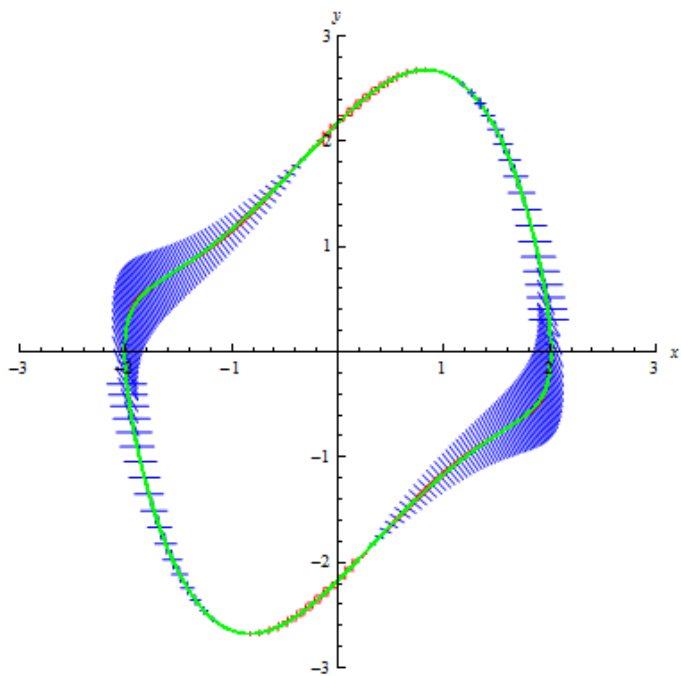

Fig. 17    A-portrait of the Van der Pol limit cycle

It is remarkable that the A-portrait of the Van der Pol limit cycle shows this limit cycle does not have "real" attractiveness everywhere. When the cycle is near to the y-axes, there is only repulsiveness distributed there.

The above four examples show that, except some extremely symmetrical cases like the example E$_3$, the uneven distribution of attractiveness exists commonly on the attractor, or on the invariant set.

In fact, this uneven distribution is very important for the bifurcation of the attractor and the invariant set.

In [11] and [12], it is shown that a stable limit cycle can separated into two stable limit cycle, and that the rotation number of a stable spatial limit closed orbit can be doubled through the bifurcation.

How can these bifurcation phenomena happen when all of these closed orbits have still the attractiveness and stability? This problem is also an important reason for the author to study the Lyapunov exponents.   But, just as shown in the last section, it has been found that none of $LE_J$, $LE_O$, $LE_Y$ and GFE can provide an exact and specific reason for these bifurcations.

Fortunately, the A-portraits of these closed orbits can show intuitively that the uneven distribution of attractiveness (especially the repulsiveness) is the main reason of these bifurcations.

Consider the following A-portraits of limit closed orbits of the Silnikov equation (10).

**E$_5$**    A-portrait in the case, $a=1, b=0.6$. In this case, the system (10) has only one stable limit cycle (Fig.2).   Figure 18 shows the uneven distribution of the attractiveness on the unique limit cycle from two viewing points.

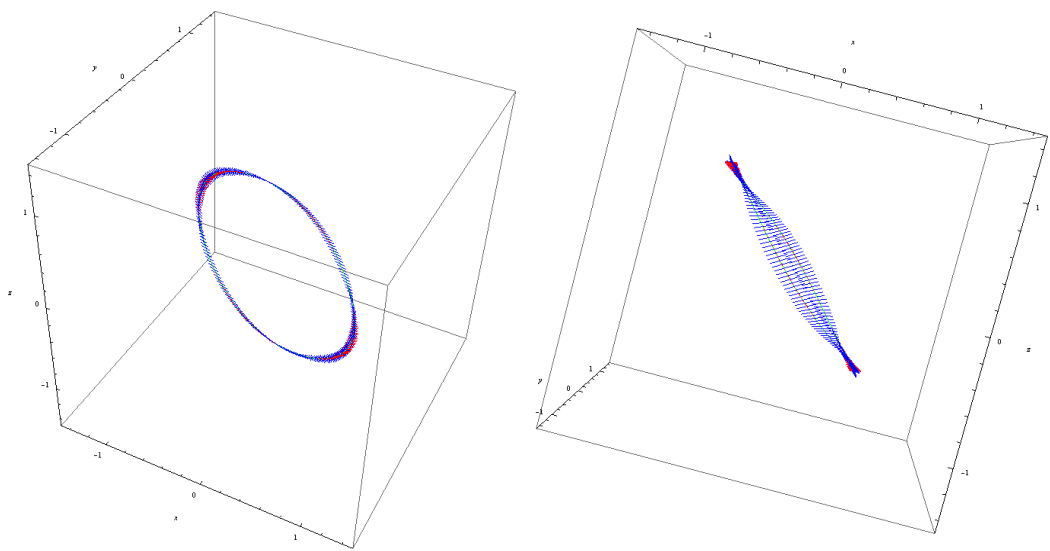

Fig. 18    A-portrait of the limit cycle in case $a=1, b=0.6$

**E$_6$**    A-portrait in the case, $a=1, b=0.4893$. In this case, the system (10) has still only one limit cycle (Fig.7), but it will be separated into two when the parameter $b$ is getting little smaller.

The A-portrait shows that the uneven of the distribution of the attractiveness is getting stronger, especially, the repulsiveness is obviously stronger (see Fig. 19).   Just the stronger repulsiveness leads the separation of the limit cycle along the repulsive direction.

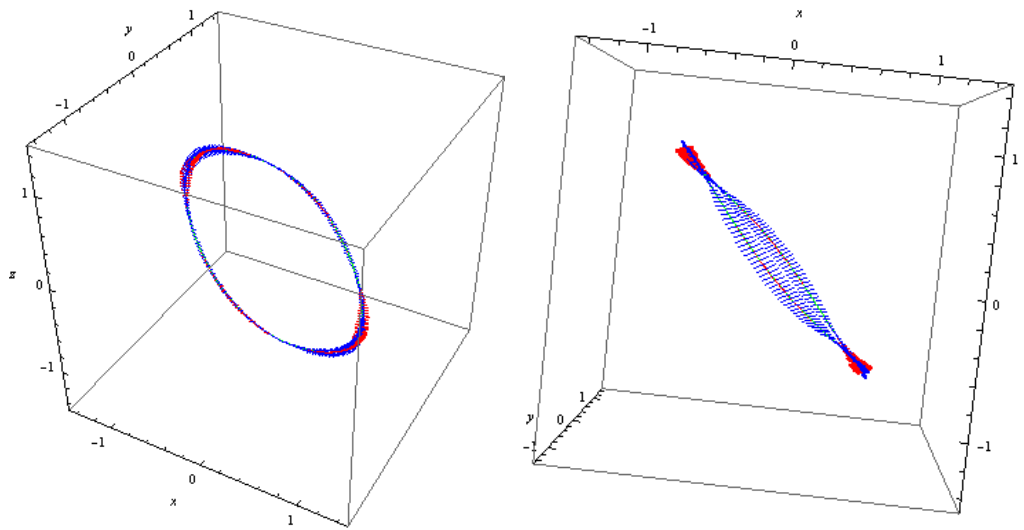

Fig.19  A-portrait of the limit cycle in case $a = 1, b = 0.4893$

**E$_7$**  A-portrait in the case, $a = 1, b = 0.4892$. In this case, the number of the limit cycles becomes two. They are slightly separated (Fig. 8). Figures 20.1 shows the uneven distribution of the attractiveness and repulsiveness on these limit cycles. Figures 20.2 is two enlarged parts of the A-portrait.   It shows the complicated relation of two limit cycles on their attractiveness and repulsiveness

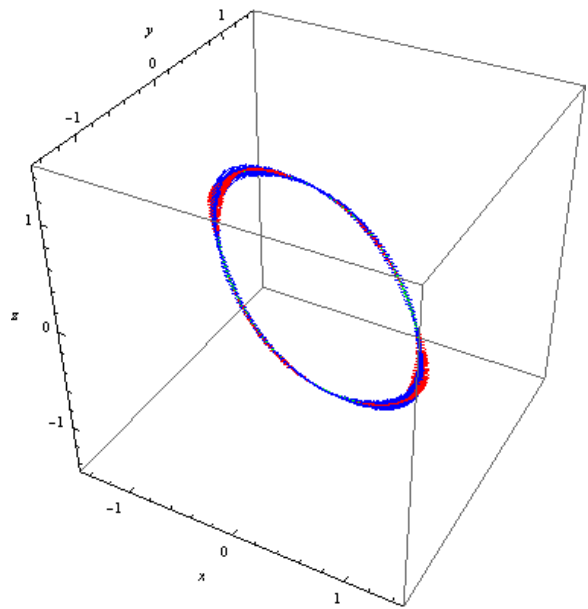

Fig. 20.1   A-portrait of double limit cycles in the case $a = 1, b = 0.4892$

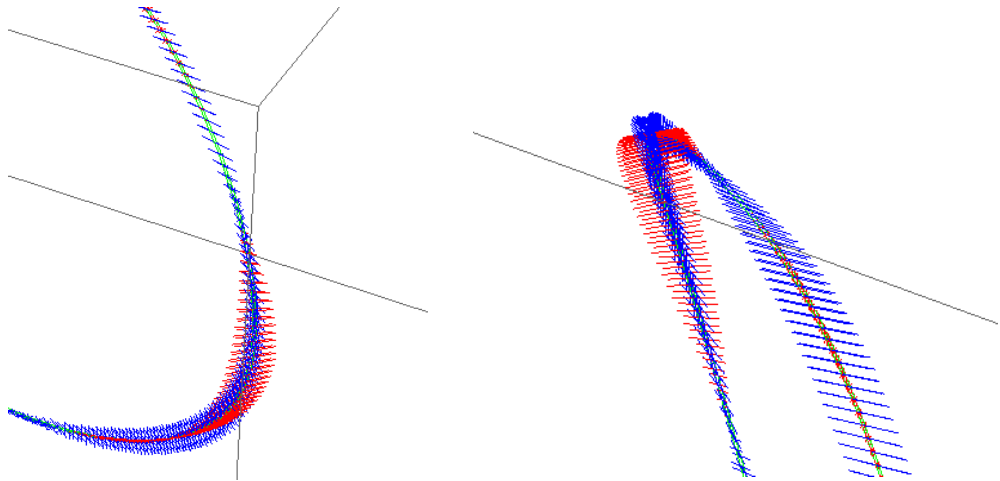

Fig. 21.2    Two enlarged parts of A-portrait

**E₈**    A-portrait in the case,  $a = 1, b = 0.4800$. In this case, two limit cycles are separated more obviously.    From two viewing points, figures 22.1 and 22.2 show their complicated relation.

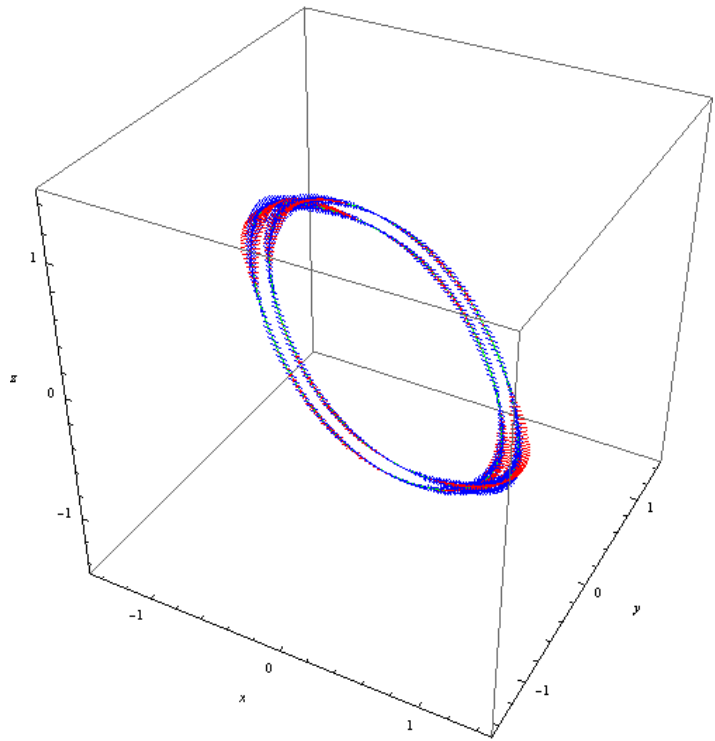

Fig.22.1    A-portrait of separated limit cycles in the case,  $a = 1, b = 0.4800$

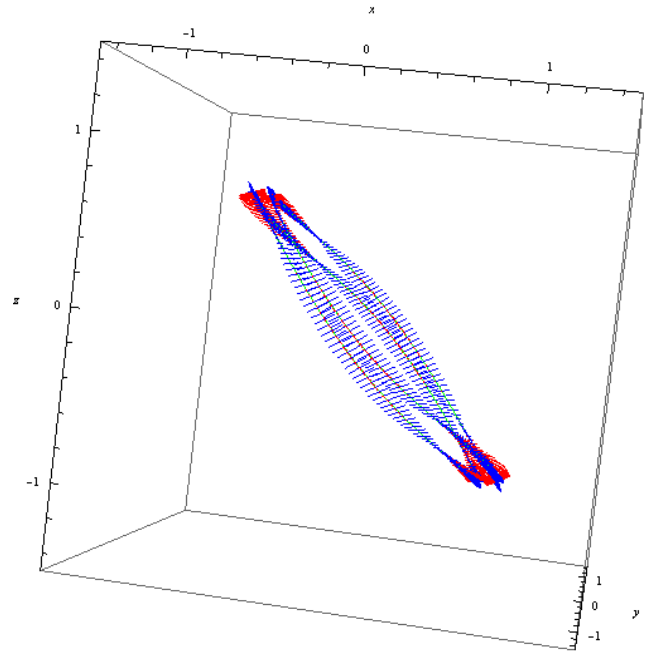

Fig.22.2  The complicated relation between two separated limit cycles

**E₉**  A-portrait in the case, $a = 1, b = 0.3995.$ When the parameter $b$ changes from 0.3996 to 0.3995, two limit cycles undergo another bifurcation, that is, the period-doubling, the rotation number of each cycle changes from one to two. Figures 23.1 shows the A-portrait of the two period-doubled closed limit orbits, and Figure 23.2 is the enlarged part of this portrait, from which the detailed complicated distribution of the attractiveness can be seen.

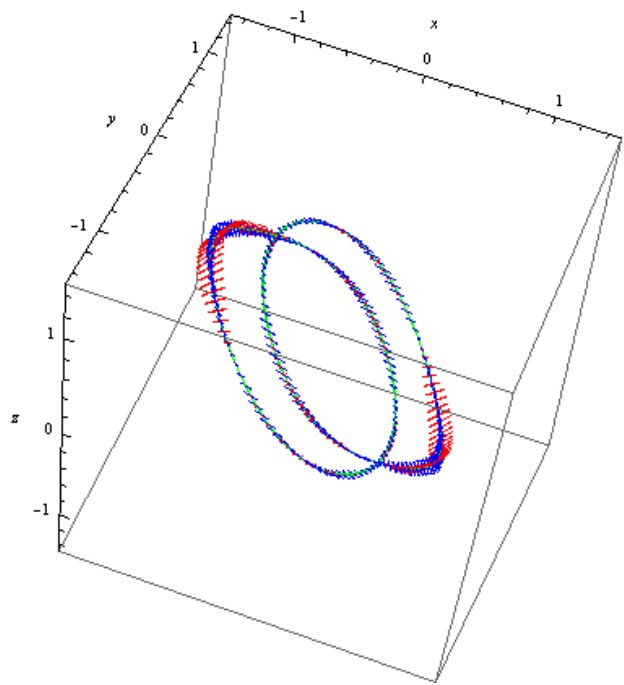

Fig.23.1  A-portrait of two period-doubled limit closed orbits

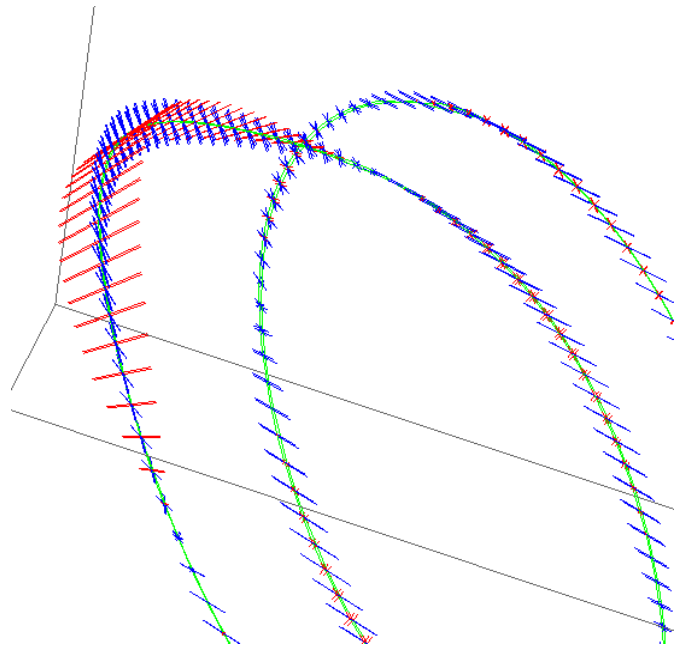

Fig. 23.2   The enlarged part of the A-portrait

**E₁₀**   A-portrait in the case, $a=1, b=0.3920$ (Fig. 24).   In this case, two period-doubled limit closed orbits can be seen clearly (ref. Fig.10).

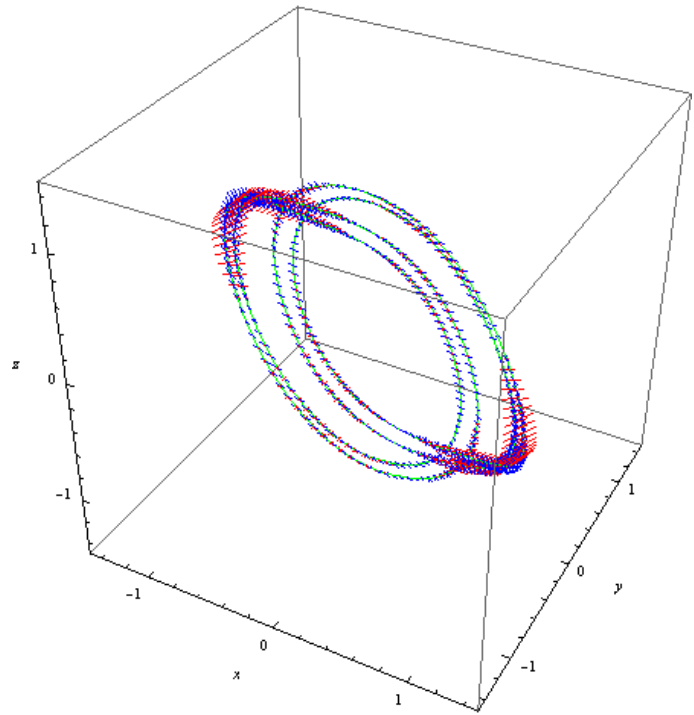

Fig. 24   A-portrait of two limit closed orbits of rotation number two

**E₁₁**   A-portrait in the case, $a=1, b=0.3338$. In this case, the attractor is a spatial limit closed orbit of rotation number 13 (Fig.11).   Its A-portrait (Figrue 25) shows the complicated distribution of the attractiveness and repulsiveness.

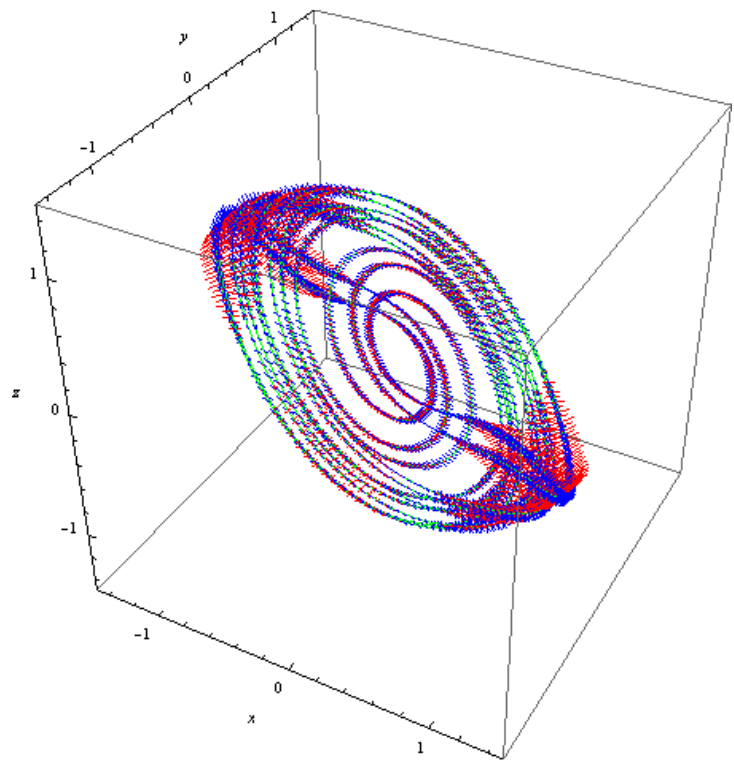

Fig. 25    A-portrait of the spatial limit closed orbit of rotation number 13.

From $E_5$ to $E_{11}$, the A-portraits for the limit closed orbits of the Silnikov equation (10) show specifically that, the uneven distribution of the attractiveness, especially the local repulsiveness, becomes stronger and stronger when the rotation number of the closed orbit is getting higher. The author thought it might explain why the leading exponent of $LE_J$ and GFE of these closed orbits with higher rotation number becomes positive.

In the following section, it will be shown that the A-portrait may provide a tool for remerging some structure hidden in a complicated attractor or in an invariant set.

## 5. Structures Hidden in an Attractor

In [12], it is shown that when the parameter $b\ (0<b<1)$ changes from large to small, the structure of the attractor of the Silnikove equation (10) may change suddenly from complicated one, such as a strange attractor, to a very simple one, such as a closed orbit with finite rotation number. Though there should be some intrinsic relation between the complicated structure and the simple structure, but the relation is usually be covered by the phase portrait of the complicated attractor, which is drawn commonly with a colored long trajectory which should have some thickness.

For instance, consider the phase portrait of the attractor in the case $a = 1, b = 0.3342$. In this case, the attractor is with an unknown complicated structure (see Fig. 26)

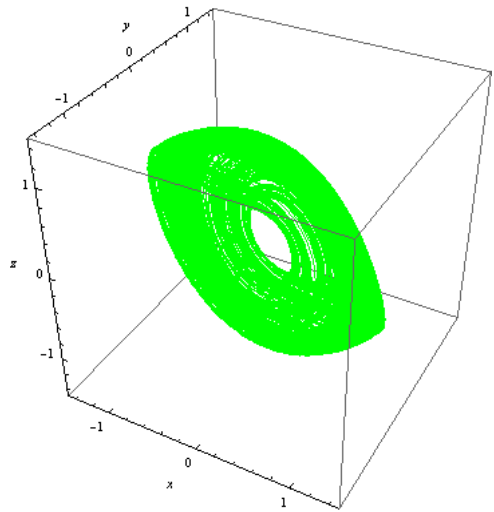

Fig. 26    The attractor with unknown complicated structure, $a = 1, b = 0.3342$

This protrait is drawn with a trajectory in the long time interval $[0, 10000]$.

When the parameter changes a little to $a = 1, b = 0.3341$, the attractor becomes a closed limit orbit with rotation number 13 (see Fig. 27)

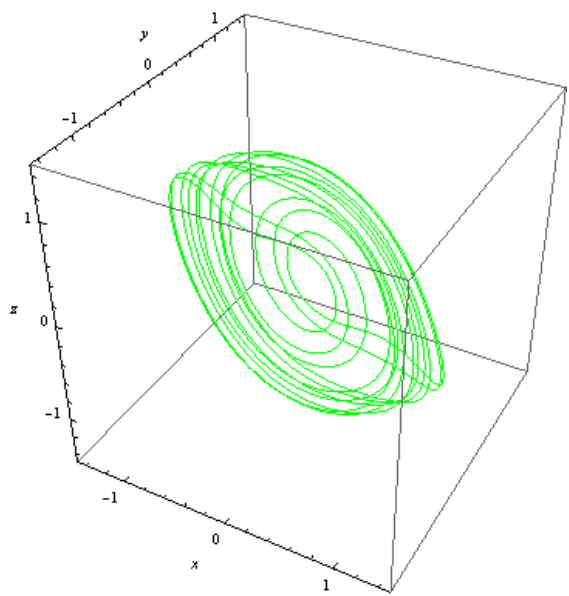

Fig.27    The closed limit orbit with rotation number 13, $a = 1, b = 0.3341$

Figures 28 and 29 are the A-portraits of the above two attractors representing in two viewing point respectively.

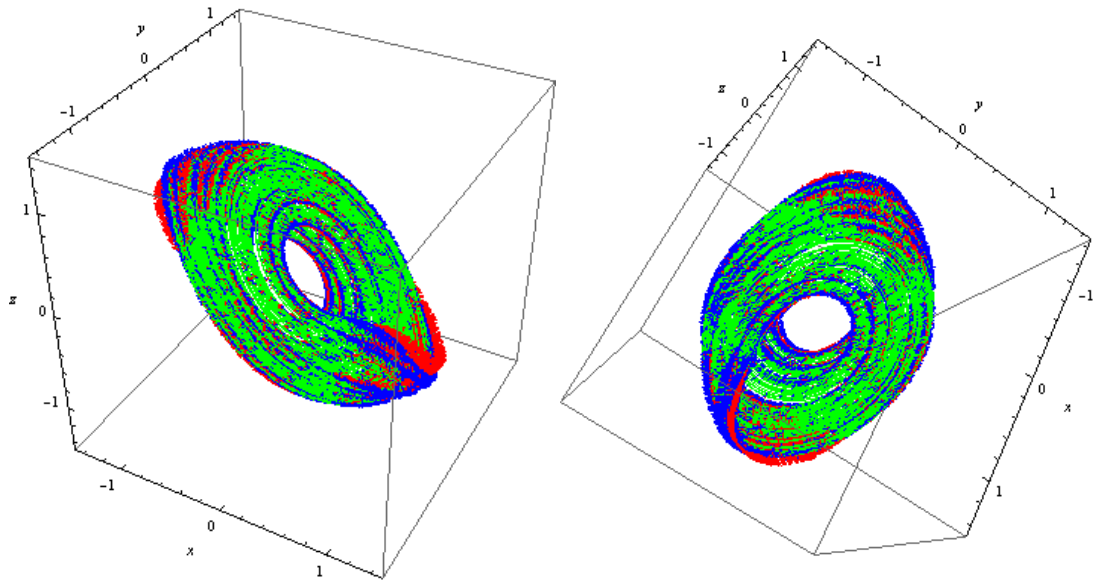

Fig.28   A-portrait of the attractor in the case   $a = 1, b = 0.3342$

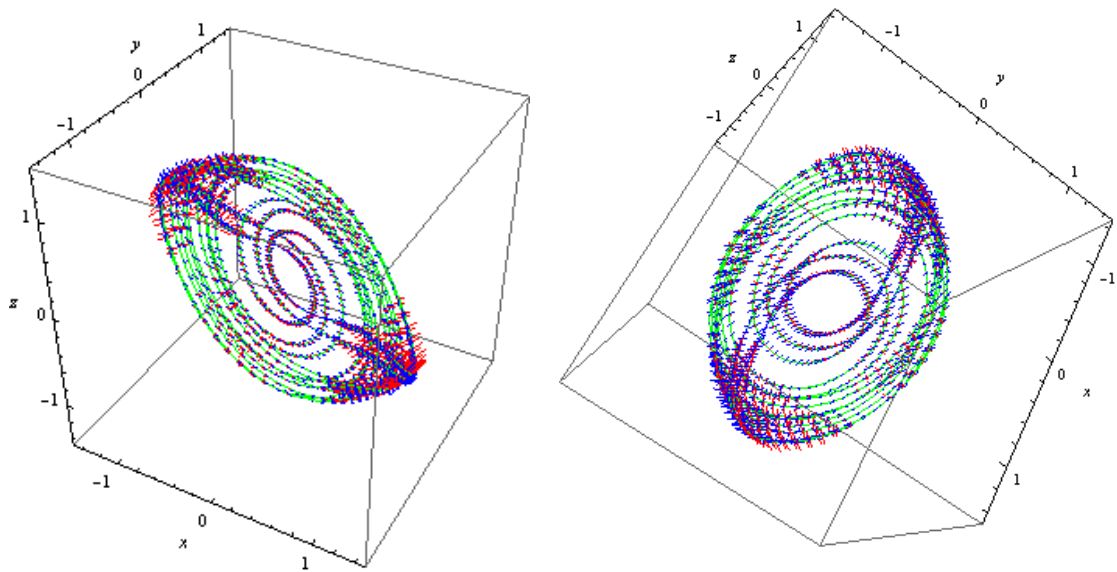

Fig.29   A-portrait of the closed orbit in the case   $a = 1, b = 0.3341$.

Comparing Fig.28 and Fig.29, it is easy to see that the structure of the A-portrait of the closed orbit of the case $a = 1, b = 0.3342$ emerges from the A-portrait of the complicated attractor of the case $a = 1, b = 0.3341$. But this structure is hidden in the commonly used phase portrait ( Fig. 26) of the complicated attractor.

Notice that $T=1$ and $m=10000$ are chosen for drawing the A-portrait of the complicated attractor, the 10000 drawing points in this portrait are evenly distributed in the long time interval $[0, 10000]$. Therefore, the emerged structure indicates that most of drawing points with stronger attractiveness or repulsiveness are centrally distributed on this structure, that is, this structure is with stronger cohesiveness in the whole attractor.

The hidden structure indicates there must be an interesting relation between the chaotic solution of an orbit on the complicated attractor and the hidden periodic solution of the closed orbit.

Figure 30 shows partially the interesting relation through the $x$ components of the two solutions in the time interval $[0,1000]$, where the horizontal axes represents time $t$, the vertical axes represents the $x$ component, the upper one is of a trajectory on the complicated attractor, and the lower one is of the periodic solution (the period is 84.1932). A phase shift is used for the comparison.

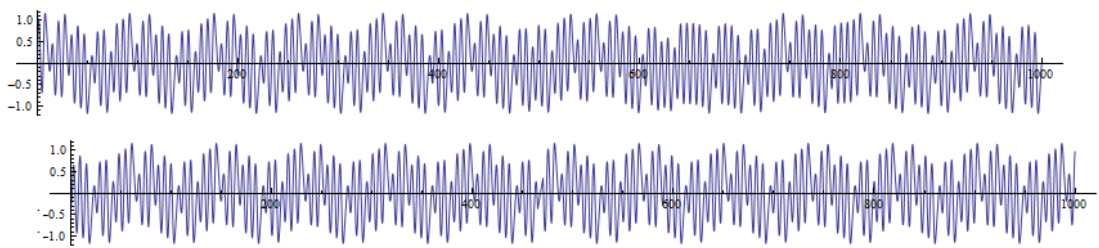

Fig.30  Comparison between the chaotic solution and the hidden periodic solution

The $y$ and $z$ components of the two solutions has the same relation.

This fact gives an answer to the problem why the complicated attractor can change into the simple closed orbit when the parameter has only a little change.

The same phenomenon happens also when the parameter changes from $a=1$, $b=0.3238$ to $a=1, b=0.3237$. When $a=1, b=0.3238$, the attractor is with a unknown complicated structure, and when $a=1, b=0.3237$, the attractor becomes a pair simple closed orbits, each of them has the same rotation number three.

Figures 31 and 32 are the A-portraits of the above two kinds of attractors representing in two viewing point respectively. The A-portraits show clearly that the simple structure of the pair of closed limit orbits in case, $a=1, b=0.3237$, are hidden in the complicated attractor of the case, $a=1, b=0.3238$.

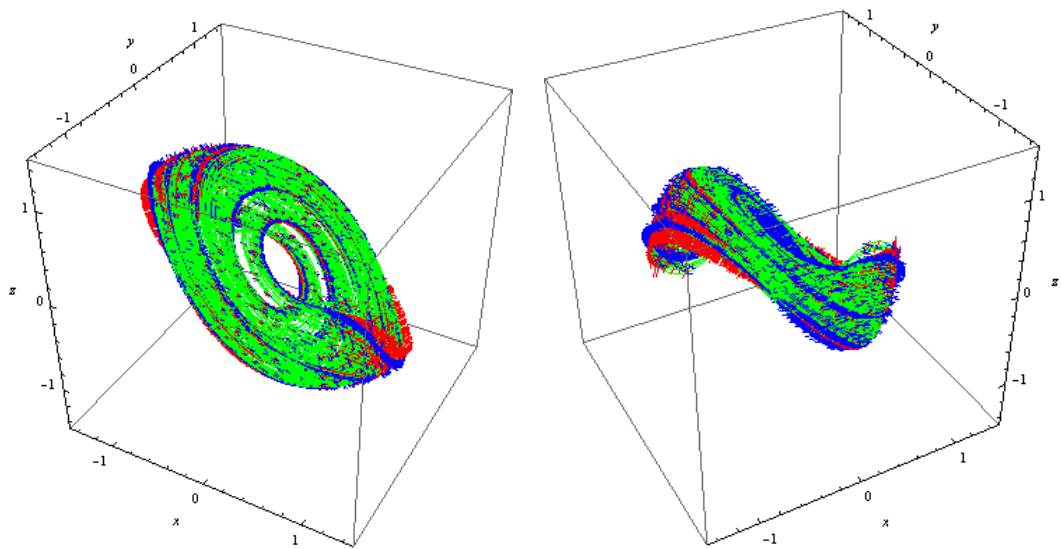

Fig.31   A-portrait of the attractor in the case   $a = 1, b = 0.3238$

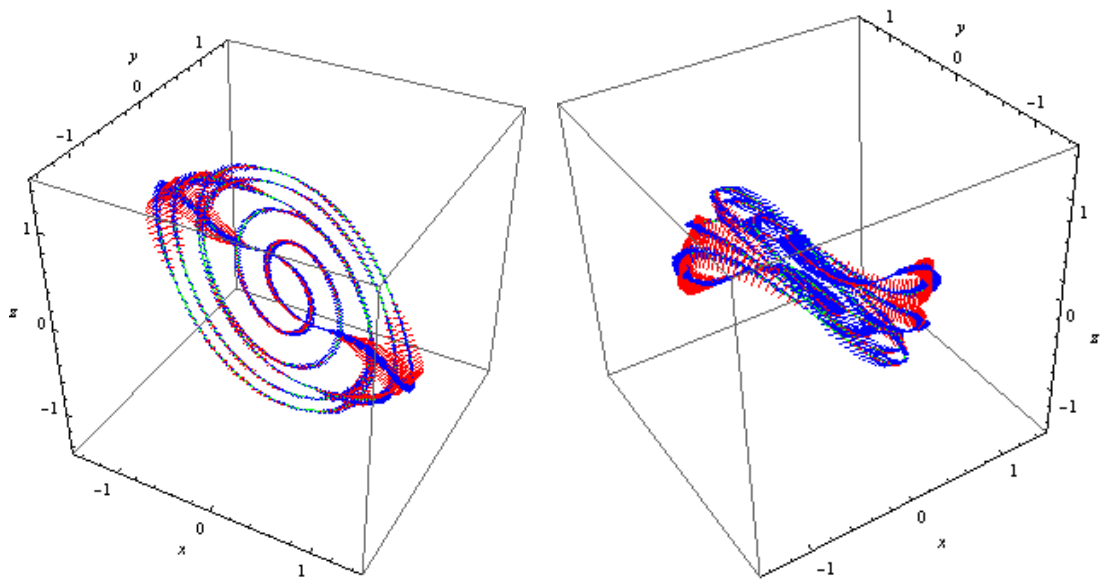

Fig.32   A-portrait of a pair of closed limit orbits in the case   $a = 1, b = 0.3237$.

It should be emphasized that the change between the closed limit orbits and complicated attractor in the above two examples does not belong to the process of the period-doubling cascade of the closed orbits.

The period-doubling cascade of the limit closed orbits is an infinite and gradual process, the structure hidden in the complicated attractor should be more interesting. This problem will be studied further.

All of above A-portrait examples are only for the attractors, none for the smallest invariant closed set. Clearly, one can study the invariant tori for Hamiltonian systems with A-portrait. But a more interesting problem related is that if it is possible that a system with dissipation can have a cluster of invariant tori with positive measurement.

Recently, William Hoover has introduced me an interesting system

$$\begin{cases} \dfrac{dq}{dt} = p \\ \dfrac{dp}{dt} = -q - \zeta p \\ \dfrac{d\zeta}{dt} = p^2 - (1 + \varepsilon \operatorname{Tanh} q) \end{cases} \quad (14)$$

which is new model of Nosé-Hoover oscillator for reconciling the time-reversible microscopic mechanics with macroscopic (irreversible) thermodynamics (ref. [13],[14],[15]). It is an improvement of another Nosé-Hoover oscillator model

$$\begin{cases} \dfrac{dq}{dt} = p \\ \dfrac{dp}{dt} = -q - \zeta p \\ \dfrac{d\zeta}{dt} = p^2 - T \end{cases} \quad (15)$$

The system (15) can be treated as an oscillatory differential equation

$$\ddot{q} + \zeta \dot{q} + q = 0$$

with friction coefficient $\zeta$ which varies with a rate $\dot{\zeta} = p^2 - T$, where $T$ is the environment temperature. This rate is corresponding to the commonly used friction coefficient which is positive if the kinetic energy of the oscillator is too large and becomes negative if the kinetic energy is too small (relative to the environment temperature $T$. The improved model (14) connects the temperature with the position variable $q$ as

$$T = 1 + \varepsilon \operatorname{Tanh} q.$$

Different to the Silnikov equation (10) and the Lorenz equation (11), which are dissipative in the whole phase space $R^3$ (when the parameters are positive) and have three equilibrium points respectively, both system (14) and system (15) are dissipative only in the half phase space where $\zeta > 0$, and both of them have no equilibrium point if $T \neq 0$.

Clearly, the fact that the system has no equilibrium state when $T \neq 0$ is an interesting property which is closer to the true nature of thermodynamics.

The above-mentioned differences increase the difficulty of the exact qualitative research for the systems (14) and (15).

Anyway, the enormous numerical results have revealed the system (14) has a series of interesting smallest invariant closed sets and attractors (ref. [13]).

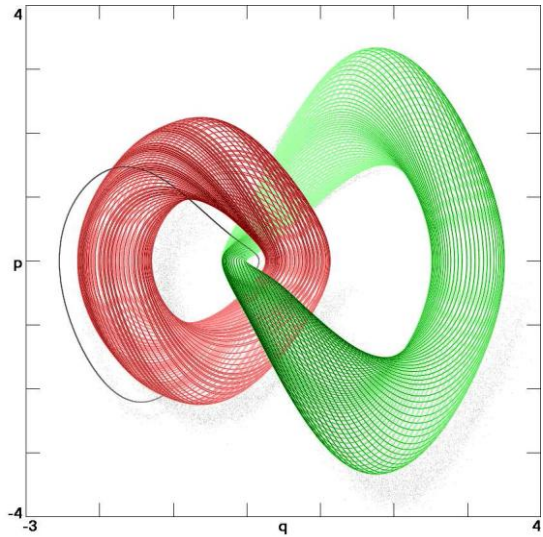

Fig.33　The q-p projection of three interlocked phase-space structures

Especially, in [13], it is found that the system (14) has three interlocked phase-space structures, two clusters of invariant tori and a dissipative limit cycle,　when $\varepsilon$ = 0.42. Figure 32 (copied from [13]) shows the q-p projection of the limit cycle and two invariant tori.

The numerical results show that both of two clusters of invariant tori may be distributed densely and may have positive measurement (see Figure 33, copied from [13]

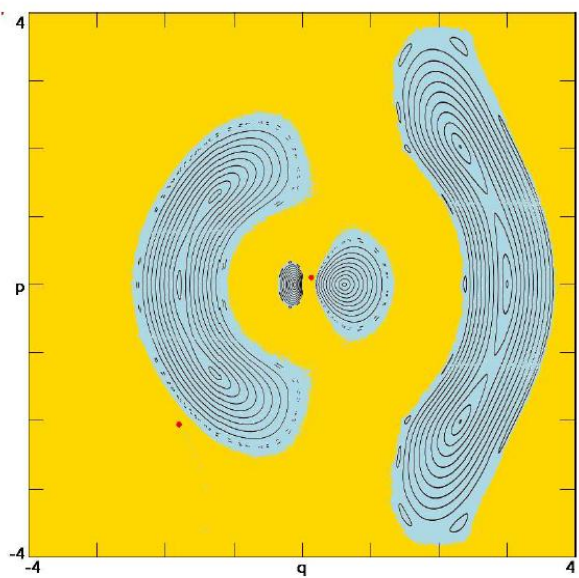

Fig.34　The detailed cross section of $\zeta = 0$ of the three structures

It is believed that just the attractiveness of the limit cycle and the positive measurements of two clusters of invariant tori makes their numerical search possible.

It is an interesting fact that the measurement of the volume enclosed by an invariant torus of the system (14) must be invariant under the flow, though this system is not conservative in the region where this volume is located. This is possible only under a very delicate balance between the divergent and convergent of the volume and between the attractiveness and repulsiveness on and inside the invariant torus.

In the case $\varepsilon = 0.42$, figure 35, 36 and 37 presents respectively the A-portrait of the limit cycle, the A-portrait of an invariant torus which is passing through the point $(-2.25, 0, 0)$, and the A-portrait of another invariant torus which is passing through the point $(2.53, 0, 0)$.

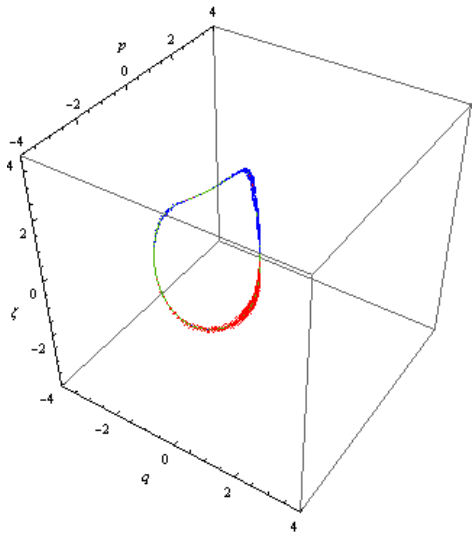

Fig. 35    A-portrait of the limit cycle

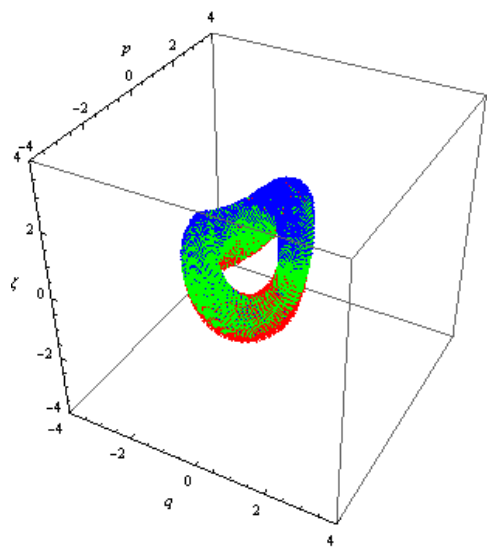

Fig. 36    A-portrait of an invariant torus passing through $(-2.25, 0, 0)$

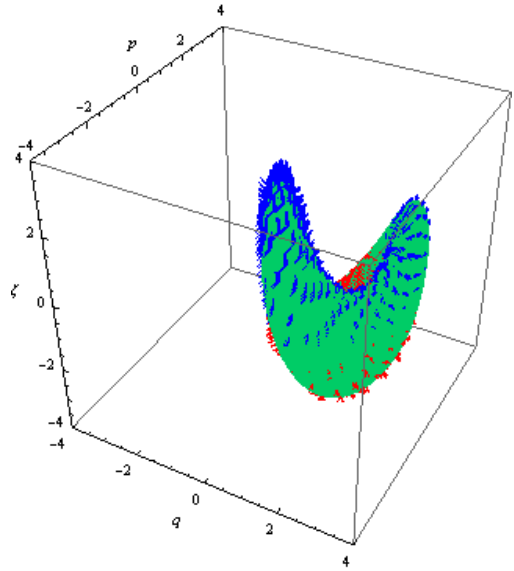

Fig. 37　A-portrait of an invariant torus passing through $(2.53, 0, 0)$

　　Putting three A-portraits together according their positions in the phase space, figures 38, 39 and 40 shows the complicate relation in attractiveness and repulsiveness between three interlocked invariant sets from three different viewing points.

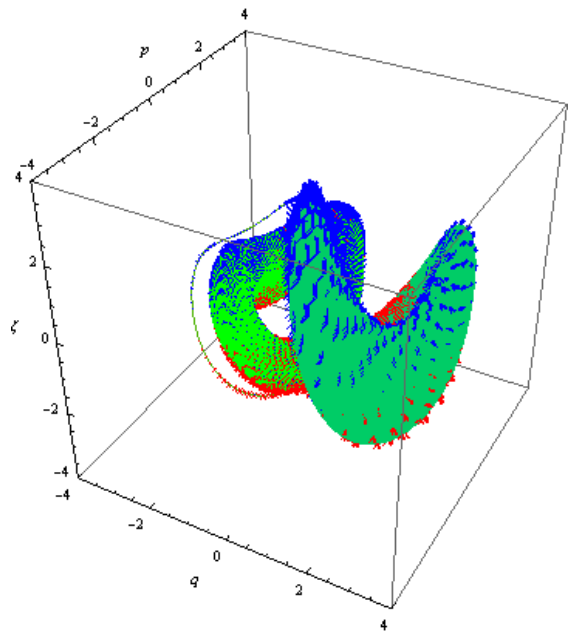

Fig.38　A-portrait of three interlocked invariant sets (1)

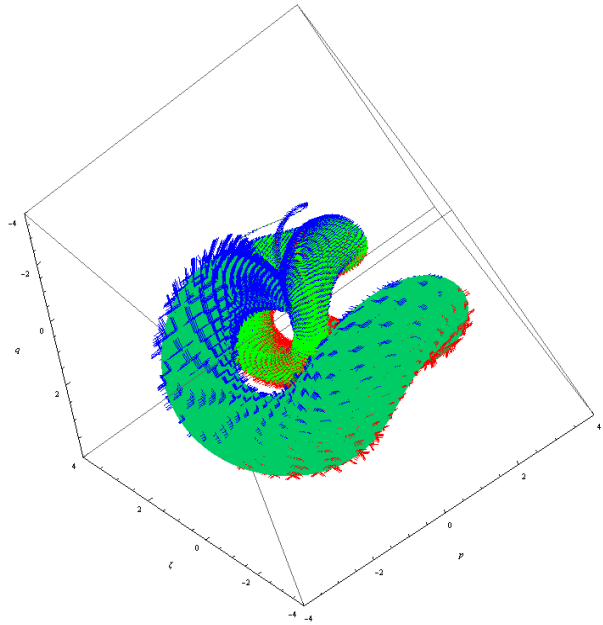

Fig.39    A-portrait of three interlocked invariant sets (2)

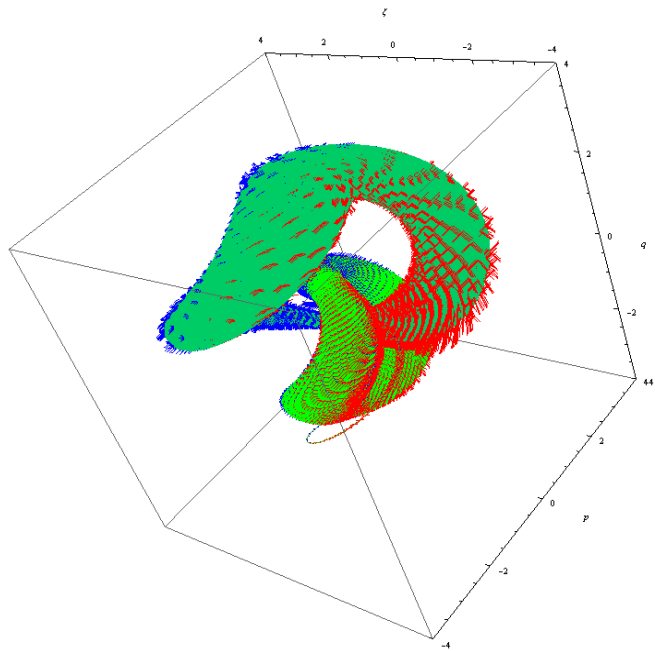

Fig.40    A-portrait of three interlocked invariant sets (3)

In drawing the A-portrait for each torus, a long trajectory corresponding long time interval $[0, 8000]$ is used, and the 8000 drawing points of attractiveness are chosen which are evenly distributed along the long time interval.    The A-portraits of two tori show that those drawing points and corresponding attractiveness are not distributed evenly on two tori.    This fact means some structures are hidden on two tori.    It is possible that these structures may imply some interesting phenomena when the parameter of the system changes.    In fact, the authors of paper [13] have studied some bifurcation phenomena of this system. They do have alobtained a series interesting results. Hope the A-portraits may provide some help for understanding their results.

# 6. Conclusion

This paper has shown that the generalized Floquet exponent may provide a more exact definition for the Lyapunov exponents for the attractor and other invariant set, though it has still some difficulties about the objectiveness.

The suggested A-portrait for an attractor or for a smallest invariant closed set may provide some more information in attractiveness of these objects for understanding the mystery process of the bifurcation phenomena of them, especially, the examples show this portrait may emerge some simple structures hidden in a complicated attractor or in an invariant set, the hidden structure may predict the possible results in a bifurcation process.

The suggested calculations are strongly relied on the development of modern computer and numerical technique.

# Acknowledgements

When the author wrote this paper, Professor William Graham Hoover happened to discuss with me about the Lyapunov exponents. Especially, he introduced me the recent research results on the Nosé-Hoover oscillator. From the discussion with him and from the joint research paper of Julien Clinton Sprott and him, I have learned a lot. Besides, he has also warmly support me to complete the present paper.

To him and to Julien Clinton Sprott who has made great contribution on their joint work, I show my hearty thanks and best regards.